\documentclass[11pt, a4paper]{amsart}

\usepackage{amsmath, amsthm, amssymb}
\usepackage{graphicx}
\usepackage{enumerate}
\usepackage[top=2.2cm, bottom=2.1cm, left=2.1cm, right=2.2cm]{geometry}

\usepackage{stmaryrd}  

\newtheorem{theorem}{Theorem}[section]
\newtheorem{corollary}[theorem]{Corollary}
\newtheorem{lemma}[theorem]{Lemma}
\newtheorem{proposition}[theorem]{Proposition}
\theoremstyle{definition}
\newtheorem{definition}[theorem]{Definition}
\newtheorem{example}[theorem]{Example}
\newtheorem{remark}[theorem]{Remark}

\numberwithin{equation}{section}

\usepackage{color}

\usepackage[T1]{fontenc}

\usepackage{algorithm, algpseudocode}

\makeatletter
\renewcommand{\p@algorithm}{\arabic{algorithm}\expandafter\@gobble}
\makeatother

\newcounter{step}[algorithm]
\newcommand\STEP[2][\(\triangleright\)]{%
	\refstepcounter{step}
	\vskip 0.25\baselineskip
	\item[]\hskip -\algorithmicindent #1 \textbf{Step \arabic{step}}%
	\ifthenelse{\equal{\unexpanded{#2}}{}}{}{ (\texttt{#2})}%
	\textbf{.}%
}
\newenvironment{algo}{\algo}{}
\def\algo#1\end{%
	\noindent\fbox{%
	\begin{minipage}[b]{\dimexpr\columnwidth-\algorithmicindent\relax}
	\begin{algorithmic}
	#1
	\end{algorithmic}
	\end{minipage}
	}%
\end}

\makeatletter
\renewcommand{\fnum@algorithm}{\fname@algorithm}
\makeatother

\usepackage{mathrsfs}

\usepackage{graphicx}
\DeclareFontEncoding{LS1}{}{}
\DeclareFontSubstitution{LS1}{stix}{m}{n}
\DeclareSymbolFont{stixletters}{LS1}{stix}{m}{it}
\DeclareMathAccent{\cev}{\mathord}{stixletters}{"91}
\DeclareMathAccent{\vec}{\mathord}{stixletters}{"92}
\DeclareMathAccent{\vecev}{\mathord}{stixletters}{"95}
\DeclareMathAccent{\cevstar}{\mathord}{stixletters}{"91}

\DeclareMathOperator*{\argmin}{argmin}

\makeatletter
\newsavebox\myboxA
\newsavebox\myboxB
\newlength\mylenA

\newcommand*\xoverline[2][0.75]{%
    \sbox{\myboxA}{$\m@th#2$}%
    \setbox\myboxB\null
    \ht\myboxB=\ht\myboxA%
    \dp\myboxB=\dp\myboxA%
    \wd\myboxB=#1\wd\myboxA
    \sbox\myboxB{$\m@th\overline{\copy\myboxB}$}
    \setlength\mylenA{\the\wd\myboxA}
    \addtolength\mylenA{-\the\wd\myboxB}%
    \ifdim\wd\myboxB<\wd\myboxA%
       \rlap{\hskip 0.5\mylenA\usebox\myboxB}{\usebox\myboxA}%
    \else
        \hskip -0.5\mylenA\rlap{\usebox\myboxA}{\hskip 0.5\mylenA\usebox\myboxB}%
    \fi}
\makeatother

\newcommand{\arcangle}{%
  \mathord{<\mspace{-9mu}\mathrel{)}\mspace{2mu}}%
  }
  \renewcommand{\angle}{\arcangle}  
  
\usepackage{tikz-cd}

\usepackage{auto-pst-pdf}
\usepackage{pst-plot}
\usepackage{pdftricks}
\begin{psinputs}
\usepackage{pstricks,pst-plot}
\end{psinputs}

\title[Alternating Bregman projections]{Alternating Bregman projections and convergence of the EM algorithm}
\author[D. Noll]{Dominikus Noll}
\address{Institut de Math\'ematiques, Universit\'e de Toulouse, France}
\email{{\tt dominikus.noll@math.univ-toulouse.fr}}

\begin{document}

\begin{abstract}
We investigate convergence of alternating Bregman projections between non-convex sets and prove convergence to a point in the intersection, or to
points realizing a gap between the two sets.  The  speed of convergence is generally  sub-linear, but may 
be linear under transversality.
We apply our analysis to prove
convergence of versions of
the expectation maximization algorithm for non-convex parameter sets.

\vspace{.1cm}
\noindent
{\bf Keywords.}
Alternating Bregman projections $\cdot$ definable sets $\cdot$ o-minimal structures 
$\cdot$ EM algorithm $\cdot$ {\it em}-algorithm $\cdot$ Kullback-Leibler divergence   $\cdot$ information geometry

\vspace{.1cm}
\noindent
{\bf AMS Classification}{
65K05 $\cdot$ 49J52 $\cdot$ 62D10 $\cdot$ 62B11}

\end{abstract}

\maketitle

\section{Introduction}
Given a finite dimensional euclidean space $\mathbb R^d$ with inner product $\langle \cdot,\cdot\rangle$ and
induced norm $\|\cdot\|$,  and a convex function $f:\mathbb R^d \to \mathbb R \cup \{+\infty\}$
of Legendre type \cite[Sect. 26]{R_legendre}, \cite{BB_legendre}, the {\it Bregman distance} or {\it Bregman divergence} associated with $f$ is defined as
\begin{equation}
\label{D}
D(x,y) = \left\{
\begin{array}{ll}
f(x) - f(y) - \langle \nabla f(y),x-y\rangle & \mbox{ if $y \in {\rm int}({\rm dom} f)=:G$ } \\
+\infty & \mbox{ otherwise}
\end{array}
\right.
\end{equation}
Given a closed subset $B \subset {\rm dom} f$, the left Bregman distance to $B$, and the left Bregman projection onto $B$,  are defined through
\begin{equation}
\label{left}
\cev{D}_B(a) = \min \{D(b',a): b'\in B\}, \quad   \cev{P}_B(a)= \argmin \{D(b',a): b'\in B\},   
\end{equation}
where the operator $\cev{P}_B$ is set-valued. Analogously, the right Bregman distance to, and projector 
onto,  a closed set $A \subset G$ are
\begin{equation}
\label{right}
\vec{D}_A(b) = \min\{D(b,a'): a' \in A\}, \quad \vec{P}_A(b) = \argmin\{D(b,a'): a'\in A\}.
\end{equation}
We consider sequences $a_n\in A$, $b_n\in B$ generated by alternating left-right Bregman projections as
$$
b_n \in \cev{P}_B(a_n), \; a_{n+1} \in \vec{P}_A(b_n), \; n\in \mathbb N,
$$
and seek conditions under which these converge to a point in the intersection, $a_n,b_n\to x^* \in A \cap B$, or in the infeasible case $A \cap B=\emptyset$, to
points $a_n\to a^*\in A$, $b_n \to b^*\in B$ minimizing the Bregman distance between $A$ and $B$. We use the notation,
$$
a_n \stackrel{l}{\longrightarrow} b_n \stackrel{r}{\longrightarrow} a_{n+1} \stackrel{l}{\longrightarrow} b_{n+1},
$$
and also the index-free form, 
\begin{equation}
\label{index_free}
a \stackrel{l}{\longrightarrow} b \stackrel{r}{\longrightarrow} a^+ \stackrel{l}{\longrightarrow} b^+,
\end{equation}
with $b\in \cev{P}_B(a)$, $a^+\in \vec{P}_A(b)$, $b^+\in \cev{P}_B(a^+)$,
referring to these as {\it building blocks} of the alternating sequence. 
Note that (\ref{index_free}) gives decrease of the distance in the sense that
\begin{equation}
\label{decrease}
D(b^+,a^+) \leq D(b,a^+) \leq D(b,a). 
\end{equation}

Alternating Bregman projections were first proposed in \cite{bregman} as iterated left projections
$
a \stackrel{l}{\longrightarrow} b \stackrel{l}{\longrightarrow} a^+ \stackrel{l}{\longrightarrow} b^+
$
between closed convex sets $A,B$. In the convex case substantial literature on Bregman projections is available,
see e.g., 
\cite{BB_legendre,BB_joint,bauschke_combettes,BCN,forward,butnariu,byrne1,censor,csiszar1,csiszar2,iusem}. 
Work addressing the non-convex case is scarce and includes for instance \cite{shawn}, even though many practical applications
use non-convex alternating Bregman procedures without proper convergence certificate. 

A strong motivation for the setup (\ref{index_free})
is that it covers instances of the {Expectation Maximization} algorithm  (EM algorithm), where the Bregman distance
specializes to the  {\it Kullback-Leibler divergence}.
This link was established in \cite{csiszar1,csiszar2}, where the authors introduce the
{\it em-algorithm},  known to coincide with the EM algorithm in the majority of cases \cite{amari}.    In \cite{csiszar1}  convergence results
for (\ref{index_free}) with $A,B$ convex are obtained. Convergence for certain non-convex instances of (\ref{index_free}) appear already in \cite{wu}. 
Here we prove
convergence for {\it definable} non-convex parameter sets, a hypothesis satisfied in practice.

Since
alternating Bregman projections include euclidean alternating projections (AP), it is worth checking
what is known in that case, as this gives an idea of what we may expect to achieve. Local convergence of non-convex AP
was proved in \cite{malick} for {\it transversal intersections}, with \cite{luke,bauschke1,bauschke2} expanding on that idea, and here
one expects linear convergence near $A \cap B$. {\it Tangential intersections} were considered
in \cite{noll, gerchberg},
and then convergence drops to sub-linear speed.  Convergence for tangential intersections is obtained under
the
{\it angle condition}, a geometric form of the {\it Kurdyka-\L ojasiewicz inequality},
controlling the mutual position of the two sets near $x^*\in A \cap B$, and this extends also to the infeasible case;  \cite{noll,gerchberg,chinese}. Presently we derive similar geometric notions for
alternating Bregman projections, including both the feasible and the infeasible case.

Along with geometry, convergence 
requires as second ingredient a weak form of {\it regularity} in at least one of the sets, or less bindingly in \cite{noll,gerchberg}, a mild form of regularity which one of the sets has
with regard to the other. Even in the euclidean case
 \cite{bauschke,DR,noll}  regularity hypotheses cannot be avoided. 
Here we rely on
the {\it three-point-inequality},  which in the non-convex setting
made its first appearance in \cite{noll}. We  obtain a suitable extension to Bregman alternating projections, and for the purpose of justification, we show  how
the three point inequality can be derived from conditions on the reach of $A,B$.

As demonstrated in \cite{noll,gerchberg} for AP, 
angle condition and regularity are truly versatile when allowed to break the symmetry between $A,B$. By that we mean that under (\ref{index_free})
regularity of $B$ with regard to $A$ does not have the same effect as regularity of $A$ with regard to $B$. Since Bregman
alternating projections are by themselves already asymmetric, 
half the theory applies to $lr$-building blocks
$a \stackrel{l}{\longrightarrow} b \stackrel{r}{\longrightarrow} a^+$, the other half
to $rl$-blocks $b \stackrel{r}{\longrightarrow} a^+ \stackrel{l}{\longrightarrow} b^+$. Fortunately,
a {\it duality principle} allows to pass from one to the other, leading to a more unified convergence theory.

We obtain
applications to the EM algorithm for discrete distributions, and for exponential families. In general only a sub-linear convergence rate
$O(k^{-\rho})$ for some $\rho\in (0,\infty)$ can be affirmed. 
Our non-smooth geometric approach has elements in common with   {\it information geometry}, a field combining Riemannian geometry, statistics and probability \cite{amari,amari_book}.

The structure of the paper is as follows. Section \ref{sect_prep} is preparatory, as is Section \ref{B_reach}, where the classical notion of reach
is extended to the Bregman setting. 
Section \ref{sect_gap} prepares the infeasible case, and in Section \ref{sect_angle_rl} the geometric angle condition  is derived from the Kurdyka-\L ojasiewicz inequality.
We also examine how the angle condition relates to tangential and transversal intersections. 
Section \ref{sect_three_point} extends the three-point-inequality from \cite{noll,gerchberg} to the Bregman setting. 
Convergence is obtained in Section \ref{sect_convergence_rl}, dual convergence in 
Section \ref{sect_dual_convergence}. The worst case speed of convergence is considered in
Section \ref{sect_speed}.
In  Section
\ref{sect_sufficient} we obtain sufficient conditions for the three point inequality from properties based on Bregman reach.
The EM algorithm is discussed in Section \ref{sect_em}. 
Terminology generally follows \cite{rock,mord}, and \cite{BB_legendre,shawn,R_legendre} are  useful for properties of Bregman distances and projectors.

\section{Preparation}
\label{sect_prep}
In this section we specify standing hypotheses, recall known results on Bregman projections, and discuss notions of reach.
Throughout we assume that $f$ is of Legendre type \cite{R_legendre,BB_legendre}, of class $C^2$ in the interior $G={\rm int}({\rm dom} f)$ of its domain, and satisfies
$\nabla^2f(x) \succ 0$ for every $x\in G$. 

\subsection{Well-posedness of the alternating sequence}
\label{sect_well}
Since $f-\langle \nabla f(a),\cdot\rangle$ is coercive for $a\in G$ by \cite[Cor. 14.2.2]{R_legendre} or \cite[Thm. 3.7]{BB_legendre}, existence 
of the left projection $\cev{P}_B(a)$ of $a\in G$
in the sense that $\emptyset \not=\cev{P}_B(a)\subset B \cap {\rm dom} f$
is assured as soon as $B \cap {\rm dom} f\not=\emptyset$ is closed in ${\rm dom} f$. This does not even require $B$ to be closed.

On the other hand, existence of the right projection $\vec{P}_A(b)$ of $b \in {\rm dom} f$ needs in the first place $A \subset G$,  but since $D(b,\cdot)$ is not coercive in general,
we assume in the alternative that $A\subset G$ is closed bounded to get $\emptyset \not= \vec{P}_A(b) \subset A$. 
With these hypotheses
alternating sequences are well-defined.

\subsection{Interiority}
\label{interior} 
The situation seems even less binding than for  iterated left projections \cite{bregman,BB_legendre}, where in order to continue left projecting from $y^+\in \cev{P}_B(y)$, one needs $y^+\in G$, a quest known as
{\it interiority}. While at first sight we do not need this here,
we have different,  yet compelling, reasons why we  also want interiority $a_k,b_k\in G$ of the alternating sequence. 

We call $B$ {\it interiority preserving} if  $B \cap {\rm dom} f\not=\emptyset$ is closed in ${\rm dom} f$ and
$\cev{P}_B(a) \subset B \cap G$ for all $a\in G$. 
It follows from
\cite[Thm. 3.12]{BB_legendre} that a closed convex  $B$ is interiority preserving iff it satisfies the constraint qualification
$B \cap G\not=\emptyset$. For non-convex $B$, the constraint qualification is necessary but no longer sufficient,
a counterexample being given  in Section \ref{sect_examples}. 

When $B$ is interiority preserving, then 
$a_k,b_k$ stay 
in $G={\rm int}({\rm dom} f)$, which has the benefit to 
enable duality, given in the next section.
However, to make full use of this, we also need
accumulation points $a^*$ of the $a_k$ and $b^*$ of the $b_k$ to stay in $G$.
Due to closedness of $A\subset G$, this is clear for the $a^*$.

Let $b^*$ be an accumulation point of the $b_k$. 
By boundedness of  $A$ we may
pass to subsequences 
$b_k \to b^*$, $a_k \to a^*$, $b_k \in \cev{P}_B(a_k)$. Suppose $b_k\to b^* \not \in {\rm dom} f$. By coercivity we have
$f(b_k)-\langle \nabla f(a^*),b_k\rangle \to \infty$, hence also
$D(b_k,a_k)= f(b_k)-f(a_k) - \langle \nabla f(a^*),b_k-a_k\rangle + \langle \nabla f(a^*)-\nabla f(a_k),b_k-a_k\rangle \to \infty$, which is absurd because by (\ref{decrease})
the sequence $D(b_k,a_k)$ is bounded. Therefore $b_k \to b^*\in {\rm dom} f$.
But then  $b^* \in \cev{P}_B(a^*)$, hence  $b^*\in G$ since $B$ is interiority preserving and
$a^*\in G$. The agreeable consequence is that the alternating sequence 
stays in a compact subset of $G$ if $A\subset G$ is closed bounded and $B$ is interiority preserving.

For any such $a_k,b_k$  we may now find a closed bounded subset $B'$ of $B$ such that $B' \subset G$ and $\cev{P}_B(a_k) \subset B'$, so that $a_k,b_k$ remain alternating
between $A$ and $B'$ with $\cev{P}_B(a_k) = \cev{P}_{B'}(a_k)$. This means the assumption
$A \subset G$ closed bounded, $B$ interiority preserving, can without loss of generality be replaced by the hypothesis $A,B \subset G$ closed bounded, which we 
adopt during the following sections.

We sketch one construction of $B'$  when $B$ is bounded. 
Find a ball $B(z,r)$ containing $A,B$ in its interior, and let $G_0 = G \cap {\rm int}B(z,r)$. Then $G_1={\rm cl}\,G_0$ is a bounded closed convex body
containing $B$, and containing $A$ in its interior. Now use a standard approximation of $G_1$ by polytopes $P_n$ from within \cite{hadwiger}, i.e.,   
$P_n \subset G_0 \subset G_1 \subset (1+\frac{1}{n})P_n$. Let $K$ be the set of all $b_k$ and all their accumulation points, then $K$  is compact, and by the above contained in $G_0$, hence
dist$(K,\partial G_0)  > 0$. Therefore some $P=P_n$ contains $K$ in its interior.  Let
$B' = B \cap P$, then $a_k,b_k$ are alternating between $A,B'$, and all accumulation points of the $b_k$ are also in $B'$.
$B'$ is convex when $B$ is, and is definable when $B$ is (see Section \ref{sect_KL}). Moreover, $B'$ is closed, because $B\cap {\rm dom} f = C \cap {\rm dom} f$ for a closed set $C$,
hence $B \cap P = B\cap {\rm dom} f \cap P= C \cap {\rm dom} f \cap P = C \cap P$ is closed.

\subsection{Legendre duality}
The conjugate $f^*$ of a function of Legendre type $f$ is again of Legendre type \cite{R_legendre,BB_legendre}, associated with $G^*= {\rm int}({\rm dom} f^*)$,  
so along with the Bregman distance $D$
generated by $f$ we also consider the distance $D^*$ generated by $f^*$ (cf. \cite[Thm. 26.5]{R_legendre}, \cite{BB_legendre}).  
For $x,y\in G = {\rm int}({\rm dom} f)$ it is known (cf. \cite[Thm. 3.7]{BB_legendre}) that
\begin{equation}
\label{fenchel}
D(x,y) = D^*(\nabla f(y),\nabla f(x)),
\end{equation}
and this gives a link between left and right projections
\begin{equation}
\label{dual_formula}
\vec{P}_A = \nabla f^* \circ \mbox{${\cev{P}}^{{}^{{}^*}}_{\nabla f(A)}$} \circ \nabla f,
\end{equation}
obtained in  \cite[Prop. 7.1]{shawn}. Here ${\vec{P}}^{{}^{{}^*}}_{\nabla f(B)}$, 
${\cev{P}}^{{}^{{}^*}}_{\nabla f(A)}$ stand for projections with regard to $f^*,D^*$. Swapping $f$ and $f^*$, we can also derive the formula
 $$\cev{P}_B = \nabla f^* \circ \mbox{${\vec{P}}^{{}^{{}^*}}_{\nabla f(B)}$}\circ \nabla f,$$ so that
$\vec{P}_A \circ \cev{P}_B = \nabla f^* \circ \mbox{${\cev{P}}^{{}^{{}^*}}_{\nabla f(A)}$} \circ \mbox{${\vec{P}}^{{}^{{}^*}}_{\nabla f(B)}$} \circ \nabla f 
= \nabla f^* \circ \mbox{${\cev{P}}^{{}^{{}^*}}_{B^*}$} \circ \mbox{${\vec{P}}^{{}^{{}^*}}_{A^*}$} \circ \nabla f$, using $\nabla f^*\circ \nabla f = id$.
Iterating this, we see that every alternating sequence (\ref{index_free}) has a mirror sequence in dual space. More precisely,
with
$$
a^* = \nabla f(b), \; b^* = \nabla f(a^+),\; a^{+*} = \nabla f(b^+), \quad A^* = \nabla f(B), \; B^* = \nabla f(A),
$$
we transform  $rl$-building blocks $b\stackrel{r}{\longrightarrow} a^+ \stackrel{l}{\longrightarrow}b^+$ between $A,B\subset G$
into $lr$-blocks $a^* \stackrel{l*}{\longrightarrow} b^* \stackrel{r*}{\longrightarrow}a^{*+}$
between $A^*,B^*\subset G^*$. Commutativity of  the following diagrams is referred to as the {\it duality principle}.

\vspace*{.1cm}
\begin{center}
\begin{tabular}{rccccccccc}
                                     & $b$ &$\stackrel{r}{\longrightarrow}$    & $a^+$ & $\stackrel{l}{\longrightarrow}$  & $b^+$  &               &  &                                    &\\
                          ${\,}^{{\,}_{\nabla\! f}} \!\!\!\!\!\!\!$           & $\big\downarrow$            &                  &  $\big\downarrow$ &                               & $\big\downarrow$ &              &            & &                                        \\
                                     & $a^*$ & $\stackrel{l*}{\longrightarrow}$ &  $b^*$ & $\stackrel{r*}{\longrightarrow}$ &$a^{*+}$&              & &                                    &   
\end{tabular}
$\quad$
\begin{tabular}{cccccccccc}
                                     & $a$ &$\stackrel{l}{\longrightarrow}$    & $b$ & $\stackrel{r}{\longrightarrow}$  & $a^+$  &               &  &                                    &\\
                            ${\,}^{{\,}_{\nabla\! f}} \!\!\!\!\!\!\!$          &                         $\big\downarrow$            &                  &  $\big\downarrow$ &                               & $\big\downarrow$ &              &            & &                                        \\
                                     & $b^*$ & $\stackrel{r*}{\longrightarrow}$ &  $a^{*+}$ & $\stackrel{l*}{\longrightarrow}$ &$b^{*+}$&              & &                                    &   
\end{tabular}
\end{center}
\vspace*{.1cm}
 It allows us to concentrate e.g. on
results for $rl$-building blocks, obtaining those for $lr$-building blocks with minor effort. This will be used repeatedly.

\subsection{Norm bounds}
\label{sect_bounds}
A consequence of non-degeneracy
$\nabla^2f(x) \succeq \epsilon I \succ 0$ of the Hessian on any compact $K\subset G$ is that we have
an estimate of the form
\begin{equation}
\label{bounds}
m^2 \|x-y\|^2 \leq D(x,y) \leq M^2 \|x-y\|^2, \; x,y\in K,
\end{equation}
with $m,M$ depending only on $f$ and $K$.  
Since Legendre functions satisfy $(\nabla f)^{-1} = \nabla f^*$, we also have an estimate of the form
\begin{equation}
\label{nabla_bounds}
l\|x-y\| \leq \|\nabla f(x)-\nabla f(y)\| \leq L\|x-y\|, \; x,y\in K,
\end{equation}
where $L$ is the Lipschitz constant of $\nabla f$ on $K$, while
$l$ is one over the Lipschitz constant of $(\nabla f)^{-1}=\nabla f^*$ on $K$.
These two imply
\begin{equation}
\label{nabla_D_bounds}
m^2L^{-2} \|\nabla f(x)-\nabla f(y)\|^2 \leq D(x,y) \leq M^2l^{-2} \|\nabla f(x)-\nabla f(y)\|^2
\end{equation}
on any such $K$. Since by our preprocessing in Section \ref{interior} we arranged $A,B \subset K \subset G$,
these norm bounds will become useful in convergence analysis.

\subsection{Uniform second-order differentiability}
\label{sect_uTJ}
The fact that $f$ is class $C^2$ on $G$ assures that it has a uniform second-order Taylor-Young expansion on every compact $K\subset G$. More precisely, given
$\epsilon > 0$, there exists $\delta > 0$ such that for all $x_0\in K$ and all $x\in G$, $\|x-x_0\| < \delta$ implies
$\left| f(x) - f(x_0) - \langle \nabla f(x_0),x-x_0\rangle - \frac{1}{2}\langle \nabla^2f(x_0)(x-x_0),x-x_0\rangle\right| \leq \epsilon \|x-x_0\|^2$. For the notion of uniform differentiability see \cite{rolewicz}.

\subsection{Kurdyka-\L ojasiewicz inequality}
\label{sect_KL}
The following definition is essential for our convergence theory.
\begin{definition}
{\rm (Kurdyka-\L ojasiewicz inequality)}.
{\rm
A lower semi-continuous function $F:\mathbb R^n \to \mathbb R \cup \{+\infty\}$ has the K\L-property at $\bar{x}\in  {\rm dom}(\partial F)$ if there exists
$\eta > 0$, a neighborhood $U$ of $\bar{x}$,  and a continuous concave function $\phi:[0,\eta) \to \mathbb R_+$, called  {\it de-singularizing function},  such that
\begin{itemize}
\item[i.] $\phi(0)=0$,
\item[ii.] $\phi$ is of class $C^1$ on $(0,\eta)$,
\item[iii.] $\phi'(s) > 0$ for $s\in (0,\eta)$,
\item[iv.] For all $x \in U \cap \{x: F(\bar{x}) < F(x) < F(\bar{x})+\eta\}$ the K\L-inequality
\begin{equation}
\label{KL}
\phi'(F(x)-F(\bar{x})) {\rm dist}(0,\partial F(x)) \geq 1
\end{equation}
is satisfied, where $\partial F$ is the limiting subdifferential. 
\end{itemize}
}
\end{definition}

\begin{remark}
We say that $F$ satisfies the 
\L ojasiewicz inequality when the de-singularizing function is $\phi'(s) = s^{-\theta}$ for some $\theta\in [\frac{1}{2},1)$, which means $\phi(s) = \frac{s^{1-\theta}}{1-\theta}$. 
\end{remark}

\begin{remark}
Information on convergence via the K\L-property  give e.g.
\cite{absil,attouch_bolte,a_b_r_s,bolte2,noll_jota,a_b_s}. 
It is well-known that definability
in an o-minimal structure \cite{dries_book,dries,dries_miller}, for short {\it definability}, implies the K\L-inequality. See \cite{kurdyka}, and for non-smooth $F$,
\cite[Thm. 11]{b_d_l_s}, where it is shown that $\phi$ may be chosen concave. 
We shall have occasion to use the small  o-minimal structure $\mathbb R_{an}$ of globally sub-analytic sets, but also large  o-minimal structures containing at least
$\mathbb R_{an,\exp}$, allowing exponential and logarithm. 
See in particular  \cite{dries_book,bierstone,coste,shiota,wilkie,miller}.
\end{remark}

\section{Bregman reach}
\label{B_reach}
This section is still preparatory, but some results are of independent interest, in particular those concerning
reach in the Bregman context, as well as concepts like Bregman geodesics encountered in information geometry  \cite{amari_book}.

\subsection{Bregman  balls}
\label{sect_rolling}
Bregman distances lead to two types of 
Bregman balls
$$
\cev{\mathcal B}(a,r) = \{x\in \mathbb R^n: D(x,a)\leq \textstyle\frac{1}{2}r^2\},
\quad \vec{\mathcal B}(b,r) = \{x\in \mathbb R^n: D(b,x) \leq \textstyle\frac{1}{2}r^2\},
$$
which generalize euclidean balls $B(x,r)$ in a natural way. Suppose $ \cev{\mathcal B}(a,r)\subset G$, and let $b\in \partial \cev{\mathcal B}(a,r)$
be a point on the boundary, then there exists a euclidean ball entirely contained in $\cev{\mathcal B}(a,r)$,
which touches the boundary $\partial \cev{\mathcal B}(a,r)$ at $b$ from within $\cev{\mathcal B}(a,r)$. We ask how the radius of this 
euclidean ball is related to $r$.

\begin{lemma}
\label{lemma1}
Let $\cev{\mathcal B}(a,r) \subset G$, and define
\begin{equation}
\label{first_kappa}
\overline{\kappa} := 
\max_{x \in \partial{\small\cev{{\mathcal B}}}(a,r)}
\frac{\lambda_{\rm max}(\nabla^2 f(x))}{\|\nabla f(x) - \nabla f(a)\|},
\quad
\underline{\kappa} := 
\min_{x \in \partial{\small\cev{{\mathcal B}}}(a,r)}
\frac{\lambda_{\rm min}(\nabla^2 f(x))}{\|\nabla f(x) - \nabla f(a)\|}.
\end{equation}
Then the euclidean ball with radius $1/\overline{\kappa}$ rolls freely inside the left Bregman ball $\cev{\mathcal B}(a,r)$. In addition, 
for every $b\in \partial \cev{\mathcal B}(a,r)$, this Bregman ball is contained in the euclidean ball with radius $1/\underline{\kappa}$
which has the same tangent hyperplane as $\cev{\mathcal B}(a,r)$ at $b$ and has its center on the same side of the tangent hyperplane as $a$.
\end{lemma}

\begin{proof}
The boundary $\partial \cev{\mathcal B}(a,r)$ of
$\cev{\mathcal B}(a,r)$ is the smooth surface 
implicitly given by the equation
$F(x) = D(x,a)-\frac{1}{2}r^2=0$.  The normal curvature of the surface
at $x\in \partial\cev{{\mathcal B}}(b,r)$
in unit tangential direction $v$ is therefore
\begin{equation}
\label{curvature}
\kappa_n(x,v) = \frac{\langle v, \nabla^2F(x)v\rangle}{\|\nabla F(x)\|} = \frac{\langle v,\nabla^2f(x)v\rangle}{\|\nabla f(x)-\nabla f(a)\|},
\end{equation}
so that $\underline{\kappa} \leq \kappa_n(x,v) \leq \overline{\kappa}$ for all $x\in  \partial{\small\cev{{\mathcal B}}}(a,r)$ and all unit $v$ from (\ref{first_kappa}).
By the Blaschke rolling theorem  \cite[§24, IV, p. 118]{blaschke}  the euclidean ball with radius
$\min_{x,\|v\|=1} 1/\kappa_n$ rolls freely inside the convex body $\cev{\mathcal B}(b,r)$, hence so does the smaller ball with radius $1/\overline{\kappa}$.

Conversely, the tangent hyperplane to $\cev{\mathcal B}(a,r)$ at $b$ being $H=\{x: \langle \nabla f(a)-\nabla f(b),x-b\rangle = 0\}$, the euclidean ball
$B(z,1/\underline{\kappa})$ with $z = b + (1/\underline{\kappa}) (\nabla f(a)-\nabla f(b))/\|\nabla f(a)-\nabla f(b)\|$ has $H$ as tangent hyperplane at $b$, and has it center
$z$ on the same side of $H$ as $a$. Since by (\ref{first_kappa}) the normal curvature $\underline{\kappa}$ of the euclidean ball is everywhere smaller than the normal curvature 
$\kappa_n(x,v)$ of the Bregman ball, we get
$\cev{\mathcal B}(a,r) \subset B(z,1/\underline{\kappa})$,  again by Blaschke's rolling theorem.
\end{proof}

On any compact subset of the interior of the domain of $f$ the eigenvalues of the Hessians $\nabla^2f(x)$  are bounded below and above by constants 
$0 < \lambda \leq \Lambda < \infty$. 
Using (\ref{bounds}) and (\ref{nabla_bounds}),
this gives
$\overline{\kappa} \leq \Lambda l^{-1} \|x-a\|^{-1}\leq {\Lambda}l^{-1} M D(x,a)^{-1/2} = {\Lambda}l^{-1}M  \sqrt{2} r^{-1}$,
and
$\underline{\kappa} \geq {\lambda} L^{-1} \|x-a\|^{-1} \geq \lambda mL^{-1} D(x,a)^{-1/2} = \lambda mL^{-1} \sqrt{2} r^{-1}$. 
We have proved the following

\begin{proposition}
\label{prop_rolling}
For every compact subset $K$ of the interior of {\rm dom}$f$ there exist constants $0 < \underline{c}\leq \overline{c}$ such that the following is true: 
If $b\in \partial \cev{\mathcal B}(a,r)$
with $\cev{\mathcal B}(a,r)\subset K$, then a euclidean ball
of radius $\underline{c}r$ is contained in $\cev{\mathcal B}(a,r)$
and touches $\partial \cev{\mathcal B}(a,r)$ at $b$ from within, and $\cev{\mathcal B}(a,r)$ is contained in a euclidean ball
with radius $\overline{c}r$ which touches $\cev{\mathcal B}(a,r)$ at $b$ from outside. 
\end{proposition}

We next consider right Bregman balls. Since these need not be convex, we need an extension of
Blaschke's rolling theorem.

\begin{lemma}
\label{walther}
Let $f$ be of class $C^{2,1}$ on $G$, and suppose $\vec{\mathcal B}(b,r) \subset G$. Let $1/r_0$
be a Lipschitz constant of $n(x) = \frac{\nabla^2f(x)(x-b)}{\|\nabla^2f(x)(x-b)\|}$ on $\partial \vec{\mathcal B}(b,r)$. 
Then a euclidean ball of radius $r_0$ rolls freely inside $\vec{\mathcal B}(b,r)$.
\end{lemma}

\begin{proof} 
Note that
$\partial \vec{\mathcal B}(b,r)$ is the $d-1$-dimensional $C^{1,1}$-sub-manifold of $\mathbb R^d$
given implicitly by the equation
$G(x) = f(b)-f(x)-\langle \nabla f(x),b-x\rangle - \frac{1}{2}r^2=0$. The outer unit normal
at a point $x$ on the boundary is therefore $n(x) = \frac{\nabla^2f(x)(x-b)}{\|\nabla^2f(x)(x-b)\|}$. Since $\partial \vec{\mathcal B}(b,r)$ is compact
and $\nabla^2f \succ 0$, $\nabla^2f(x)$ is bounded and bounded away from 0
on the boundary. Moreover,
$b\not\in \partial \vec{\mathcal B}(b,r)$, hence the denominator $\|\nabla^2f(x)(x-b)\|$ stays bounded away from 0. Finally, since by hypothesis  $\nabla^2f$ is locally Lipschitz,
the unit normal $n(x)$ is Lipschitz on the compact $\partial \vec{\mathcal B}(b,r)$. Assuming $1/r_0$ is a Lipschitz constant, it follows with the extension
\cite[Thm. 1 (v)]{walther} of Blaschke's rolling theorem
that a ball of radius $r_0$ rolls freely inside $\vec{\mathcal B}(b,r)$. 
\end{proof}

Using (\ref{bounds}) and $0 < \lambda \leq \lambda_i \leq  \Lambda < \infty$
for the eigenvalues $\lambda_i$  of $\nabla^2f(x)$ on $\partial \vec{\mathcal B}(b,r)$,
we may  again show that $r_0$ is proportional to $r$ in the rough sense of Proposition \ref{prop_rolling}. We leave the details to the reader.

\subsection{Bregman geodesics}
Let $b^+ \in \cev{P}_B(a^+)$ with $a^+\not\in B$, and define 
$
a_\lambda = (\nabla f)^{-1}(\lambda \nabla f(a^+)+(1-\lambda) \nabla f(b^+)).
$
 Then
$\cev{P}_B(a_\lambda) = b^+$ for $0 \leq \lambda < 1$; cf. \cite[Prop. 3.2]{shawn}. We call the curve $a_\lambda$ the left Bregman perpendicular to $B$ at $b^+$ 
from $a^+\in G\setminus B$. 
The direction
$d=\nabla^2 f(b^+)^{-1} (\nabla f(a^+)-\nabla f(b^+))$ is tangent to the curve $a_\lambda$ 
at $a_\lambda|_{\lambda=0}=b^+$; see also \cite{shawn}. In addition,   
$d$ is  normal to the tangent hyperplane $H=\{b: \langle \nabla f(a^+)-\nabla f(b^+),b-b^+\rangle = 0\}$ to $B$ at $b^+$ in the euclidean geometry
$\|x\|_{b^+}^2 = \langle x, \nabla^2 f(b^+)x\rangle = \langle x,x\rangle_{b^+}$. The curve $a_\lambda$ consists of 
those points which satisfy $\cev{P}_H(a_\lambda) = b^+$ and are
on the same side of $H$ as $a^+$.  Perpendiculars $a_\lambda$ are also known as left Bregman geodesics.
Note that a  non-smooth $b^+\in B$ may be left projected on from
different points lying on different geodesics. Those may then be distinguished by their tangents $d$ at $b^+$.

We next investigate right Bregman geodesics. As we shall see, this may be based on the dual
formula (\ref{dual_formula}).
\begin{lemma}
\label{lemma3}
Let $a^+ \in \vec{P}_A(b)$ with $b\not\in A$, and define $b_\lambda = \lambda b + (1-\lambda) a^+$. Then $\vec{P}_A(b_\lambda) = a^+$ for
$0 \leq \lambda < 1$.
\end{lemma}

\begin{proof}
We put $a^{*+} = \nabla f(b)$, so that
$b = \nabla f^*(a^{*+})$, also $b^{*+} = \nabla f(a^+)$ and $B^*=\nabla f(A)$.
Then $b^{*+} = \nabla f(a^+) \in \nabla f \circ \vec{P}_A(\nabla f^*(a^{*+}))
= 
\mbox{${\cev{P}}^{{}^{{}^*}}_{\nabla f(A)}$}(a^{*+})  =\mbox{${\cev{P}}^{{}^{{}^*}}_{B^*}$}(a^{*+})$ 
by (\ref{dual_formula}). 
This means we are in the situation of the left Bregman projection  above, with $f$ replaced by $f^*$ 
and $B$ replaced by $B^*$.
The left Bregman perpendicular to $B^*$ at $b^{*+}$, being $a_\lambda^* = \nabla f(\lambda \nabla f^*(a^{*+}) + (1-\lambda)\nabla f^*(b^{*+}))$,
is defined for $\lambda \in [0,1]$, and for $0 \leq \lambda < 1$, $\mbox{${\cev{P}}^{{}^{{}^*}}_{B^*}$}(a^{*}_\lambda)= b^{*+}$ is single-valued.
Reading the dual formula (\ref{dual_formula}) backward means $\vec{P}_A(b_\lambda) = a^+$  single valued for $0 \leq \lambda < 1$.
\end{proof}

We call $b_\lambda$  the right Bregman perpendicular to $A$ at $a^+\in \vec{P}_A(b)$ from $b\not\in A$, or right geodesic.
Right geodesics are straight lines, while left geodesics
are  curved. Let $D(b,a^+)=\frac{1}{2}r^2$, then the tangent hyperplane to $\vec{\mathcal B}(b,r)$ at $a^+\in \partial \vec{\mathcal B}(b,r)$
is $H =\{z: \langle \nabla^2f(a^+)(b-a^+),z-a^+\rangle = 0\}$. The right Bregman perpendicular to $A$ at $a^+$ in direction $d=b-a^+$
includes, and for convex right Bregman balls consists of,   those points $b_\lambda$ which satisfy $\vec{P}_H(b_\lambda) = a^+$ and lie on the same side of $H$ as $b$. The geodesic is also
the normal to $H$ in the euclidean geometry  induced by the scalar product $\langle x,\nabla^2f(a^+)y\rangle$.

\begin{remark}
In the case of the Kullback-Leibler divergence, left and right perpendiculars  are known
as $m$- and $e$-geodesics; cf. \cite{amari,amari_book}.
\end{remark}

\subsection{Normal cones}
\label{sect_left_project}
We recall the following properties of the left and right Bregman projectors.

\begin{lemma}
\label{prox_normals}
Let $a \stackrel{l}{\longrightarrow} b \stackrel{r}{\longrightarrow} a^+ \stackrel{l}{\longrightarrow} b^+$ 
be a building block. 
Then
\begin{itemize}
\item[(i)] $\nabla f(a^+)-\nabla f(b^+) \in N_B^p(b^+)$. 
\item[(ii)] $\nabla^2f(a^+)(b-a^+)\in \widehat{N}_A(a^+)$.
\item[(iii)] When $f$ is of class $C^{2,1}$, then $\nabla^2f(a^+)(b-a^+)\in N^p_A(a^+)$.
\end{itemize} 
\end{lemma}

\begin{proof}
For (i) see \cite[Prop. 3.3]{shawn}. 
This can also be derived directly from Lemma \ref{lemma1}.  Place a euclidean ball $B(z,\underline{c}r)$ such that it is contained in $\cev{\mathcal B}(a^+,r)$ and touches at $b^+$ from within, where
$D(b^+,a^+)=\frac{1}{2}r^2$. Then $b^+ \in P_B(z)$ with 
$z\in b^+ + \mathbb R_+(\nabla f(a^+)-\nabla f(b^+))$.

Concerning (ii), from $a^+\in \vec{P}_A(b)$ we derive
$a^+\in \partial \vec{\mathcal B}(b,r)$ with $\frac{1}{2}r^2 = D(b,a^+)$. The tangent hyperplane to the regular surface $\partial \vec{\mathcal B}(b,r)$ at $a^+$ has normal 
$n(a^+)=\nabla^2f(a^+)(b-a^+)$. This means $n(a^+)$ is a Fr\'echet normal to the sublevel set
$S=\{x: -D(b,x) \leq -\frac{1}{2}r^2\}$ at $a^+$, and since $a^+\in A \subset S$, $n(a^+)$ is also a Fr\'echet normal to $A$ at $a^+$.

Part (ii) may also be derived via duality. The building block $a \stackrel{l}{\longrightarrow} b \stackrel{r}{\longrightarrow} a^+$
is the image under $\nabla f^*$ of a building block $b^* \stackrel{r*}{\longrightarrow} a^{*+} \stackrel{l*}{\longrightarrow} b^{*+}$. 
For the latter $\nabla f^*(a^{*+})-\nabla f^*(b^{*+}) \in N_{B^*}^p(b^{*+})$ by part (i). 
Using $N_{B^*}^p(b^{*+}) \subset \widehat{N}_{B^*}(b^{*+})$ and the chain rule \cite[Thm. 10.6, Ex. 6.7]{rock} gives $\nabla^2f(a^+)(b-a^+) \in \widehat{N}_A(a^+)$.

Finally, concerning (iii), when $f$ is class $C^{2,1}$,  we may by Lemma \ref{walther}
place a euclidean ball $B(z,r_0)$ such that it has $a^+$ on its boundary and is contained in $\vec{\mathcal B}(b,r)$, where $D(b,a^+)=\frac{1}{2}r^2$.
Since $B(z,r_0)$ and $\vec{\mathcal B}(b,r)$ share the tangent hyperplane
$H=\{x: \langle \nabla^2f(a^+)(b-a^+),x-a^+\rangle = 0\}$ at $a^+$, the center $z$ must lie on the normal 
$a^+ + \mathbb R_+d$, $d= \nabla^2f(a^+)(b-a^+)$, which is therefore proximal.
\hfill $\square$
\end{proof}

\subsection{Bregman reach}
\label{sect_reach}
A point $b^+\in B$ is {\it projected on} if there exists $c \not\in B$ with $b^+\in P_B(c)$. Then $d:= c-b^+\in N_B^p(b^+) \not=\{0\}$,
and the reach $R(b^+,d)$ of $B$ at $b^+$ in direction $d$  is the radius $r$ of the largest ball $B(c,r)$ with center 
$c = b^++rd/\|d\|$ having no point of $B$ in  its interior (cf. \cite{federer}). It is possible that $R(b^+,d)=\infty$, when the largest ball becomes a half-space. 

We extend this classical notion of reach \cite{federer} to the Bregman setting. 
The results in Section \ref{sect_rolling} show that if $b^+$ is projected on, 
then it is also {\it left Bregman projected on}, i.e., $b^+\in \cev{P}_B(a^+)$ for some $a^+\not\in B$.

\begin{definition}
{\rm
Let $b^+\in \cev{P}_B(a^+)$ be left projected on from some $a^+\in G\setminus B$, and let $a_\lambda$ be the associated left Bregman geodesic.
The left Bregman reach $\cev{R}(b^+,d)$ of $B$  at $b^+\in B$ in direction $d=\nabla^2f(b^+)^{-1}(\nabla f(a^+)-\nabla f(b^+))$ is the largest left Bregman radius $r_\lambda$ for which
the Bregman ball $\cev{\mathcal B}(a_\lambda,r_\lambda)$ with $b^+\in \partial \cev{\mathcal B}(b_\lambda,r_\lambda)$ has no points of $B$ in its interior.
}
\end{definition} 

\begin{remark}
As in the euclidean case
the Bregman reach may be infinite, which is when  $\bigcup_{\lambda > 0}\cev{\mathcal B}(a_\lambda,r_\lambda)$ contains no point of $B$ except $b^+$.
The definition reproduces the classical definition of reach
in the euclidean case.
\end{remark}

\begin{definition}
\label{def_left_reach}
{\bf (Left Bregman reach)}.
{\rm 
The set $B$ has left Bregman reach at least $r > 0$ at $b^*\in B$, written $\cev{R}(b^*)\geq r$,  if there exists a neighborhood $U$ of $b^*$ such that
for every point $b^+ \in \cev{P}_B(a^+) \cap U$ left projected on from some $a^+\in G \setminus B$,  we have $\cev{R}(b^+,d) \geq r$ for the associated
$d=\nabla^2f(b^+)^{-1}(\nabla f(a^+)-\nabla f(b^+))$.
}
\end{definition}

\begin{definition}
{\bf (Right Bregman reach)}.
{\rm
The right Bregman reach $\vec{R}(a^+,d)$ of $A$ at $a^+\in A$ in direction $d=b-a^+\not=0$  with $a^+\in \vec{P}_A(b)$,  is the radius $r_\lambda$
of the largest right Bregman ball $\vec{\mathcal B}(b_\lambda,r_\lambda)$ with $a^+\in \partial \vec{\mathcal B}(b_\lambda,r_\lambda)$ having no points of
$A$ in its interior. $A$ has right Bregman reach at least $r > 0$ at $a^*\in A$, written $\vec{R}(a^*) \geq r$,  if $\vec{R}(a^+,d) \geq r$ on a neighborhood $U$ of $a^*$.
}
\end{definition}

\begin{remark}
Suppose $D(b,a_\lambda) = D(b^+,a_\lambda)=\frac{1}{2}r_\lambda^2$, then $D^*(b^*_\lambda,a^*) = D^*(b^*_\lambda,a^{*+})$ by duality,
as can be seen from $D^*(\nabla f(a_\lambda),\nabla f(b)) = D^*(\nabla f(a_\lambda),\nabla f(b^+)) = \frac{1}{2}r_\lambda^2$. This shows
$\cev{R}(b^+,d) = \vec{R}(a^{*+},d_*)$  with $d=\nabla^2f(b^+)^{-1}(\nabla f(a^+)-\nabla f(b^+))$ and $d_* = b^*-a^{*+}$. 
In other words,
Bregman reach is amenable to duality. 
\end{remark}

\begin{remark}
As a consequence of the results in Section \ref{sect_rolling}, we now see that if $K \subset G$ is compact, then any set $B \subset K$ with left Bregman reach
$r>0$  has classical reach at least $\underline{c}r>0$ and at most $\overline{c}r$, where $\underline{c},\overline{c}$ depends only on $K$, and inversely, if $B$ has classical reach $r$, then
it has left Bregman reach at most $\underline{c}^{-1} r$, and at least $\overline{c}^{-1}r$. Right Bregman reach and classical reach are also proportional in this rough sense when $f$ is class $C^{2,1}$. 
\end{remark}

\subsection{$1$-coercivity of $f$}
Consider a left perpendicular $a_\lambda$ at $b^+\in \cev{P}_B(a^+)$ from $a^+\not\in B$. Let $H = \{x: \langle \nabla f(b^+)-\nabla f(a^+),x-b^+\rangle=0\}$
be the tangent hyperplane to $\partial \cev{\mathcal B}(a^+,r)$ at $b^+$, where $D(b^+,a^+)=\frac{1}{2}r^2$,  $H_+$ the half space containing $a^+$, $H_{++}$ the open half space.
Suppose $a_\lambda$ is defined for
$0 \leq \lambda < \lambda_\infty$, where $\lambda_\infty \in (1,\infty]$, $D(b^+,a_\lambda)=\frac{1}{2}r_\lambda^2$. 
We ask whether the balls $\cev{\mathcal B}(a_\lambda,r_\lambda)$,  $0 \leq \lambda < \lambda_\infty$,
fill the open half space $H_{++} \cap {\rm int}({\rm  dom}f)$. 

\begin{proposition}
\label{prop2}
The following are equivalent:
\begin{enumerate}
\item[\rm (i)]
For some nonempty compact $B\subset G={\rm  int}({\rm dom} f)$, and every $a^+\not \in B$, $b^+\in \cev{P}_B(a^+)$ with perpendicular $a_\lambda$,
the balls $\cev{\mathcal B}(a_\lambda,r_\lambda)$,  $0 \leq \lambda < \lambda_\infty$
fill the half space  $H_{++} \cap {\rm int}({\rm  dom}f)$.
\item[\rm (ii)] $f$ is $1$-coercive. 
\end{enumerate}
When these are satisfied, then {\rm (i)} holds for {\rm every} such $B$, and we have $\lambda_\infty = \infty$ in {\rm every} perpendicular curve $a_\lambda$.
\end{proposition}

\begin{proof}
1)
Assume first that $f$ is $1$-coercive. Then by \cite[Prop. 2.16]{BB_legendre} $f^*$ is defined everywhere, hence so is $\nabla f^*$, hence 
$a_\lambda$ is defined for all $\lambda \geq 0$. Now let $b\in H_{++}\cap {\rm int}({\rm dom} f)$. That means $\langle \nabla f(a^+)-\nabla f(b^+),b-b^+\rangle > 0$.
By the cosine theorem for Bregman distances and the definition of $a_\lambda$ we have
$$
D(b,a_\lambda )= D(b,b^+) + D(b^+,a_\lambda) - \lambda \langle \nabla f(a^+)-\nabla f(b^+),b-b^+\rangle
$$
hence
\begin{align*}
\frac{D(b,a_\lambda)-D(b^+,a_\lambda)}{\lambda} &= \frac{D(b,b^+)}{\lambda} -\langle \nabla f(a^+)-\nabla f(b^+),b-b^+\rangle\\
&\to -\langle \nabla f(a^+)-\nabla f(b^+),b-b^+\rangle < 0
\end{align*}
as $\lambda \to \infty$, 
so that $D(b,a_\lambda) < D(b^+,a_\lambda)$ for $\lambda$ large enough,
which gives $b \in \cev{\mathcal B}(a_\lambda,r_\lambda)$. 

2)
Conversely, suppose $H_{++}\cap {\rm int}({\rm dom} f) \subset \bigcup \{\cev{\mathcal B}(a_\lambda,r_\lambda): 0 \leq \lambda < \lambda_\infty\}$. The normal curvature
of $\partial \cev{\mathcal B}(a_\lambda,r_\lambda)$ at $b^+\in B$ in unit tangential direction $v$ being
$$
\kappa_n(\lambda)=\frac{\langle v, \nabla^2f(b^+)v\rangle}{\| \nabla f(b^+)-\nabla f(a_\lambda)\| },
$$
we see that
$\kappa_n(\lambda)$ must tend to 0 as $\lambda \to \lambda_\infty$.  Indeed, by convexity left Bregman balls touching $H$ at $b^+$
have to get arbitrarily flat at $b^+$ as $\lambda$ increases  in order to contain points $b\in H_{++}$ a fixed distance away from $b^+$, while arbitrarily close to $H$. But that means 
the denominator $\|\nabla f(b^+)-\nabla f(a_\lambda)\|$ must tend to infinity, since the numerator is bounded below by $\lambda_{\min}(\nabla^2f(b^+))>0$. 
That forces $\|\nabla f(a_\lambda)\|\to \infty$ as $\lambda \to \lambda_\infty$ for every perpendicular $a_\lambda$ to $B$ at any $b^+\in B$ left projected on from some
$a^+\not\in  B$.  

Now suppose that contrary to what is claimed 1-coercivity fails.  By  \cite[Lemma 26.7]{R_legendre} this means 
 there exist $a_k$ with $\|a_k\|\to \infty$ such that
$\|\nabla f(a_k)\|\leq K < \infty$. 
Using compactness of $B \subset G$, find $b_k\in \cev{P}_B(a_k)$, and let $a_{k,\lambda}$ be the left perpendicular to $B$ at
$b_k$ from $a_k$.
Then we can find
points $\bar{a}_k$ on the perpendicular,  having $b_k = \cev{P}_{B}(\bar{a}_k)$, such that $\|\bar{a}_k\| \leq K' < \infty$, and at the same time
$\|\bar{a}_k-b_k\|\geq \epsilon > 0$ for all $k$ (the latter using $\|a_k-b_k\|\to \infty$, allowing to choose $\bar{a}_k$ in between bounded while away from $b_k$).  
Passing to subsequences, we may assume $b_k \to b^*$, $\bar{a}_k \to a^*$ with $b^*\in \cev{P}_B(a^*)$ and $\|a^*-b^*\|\geq \epsilon$, while
$\nabla f(a_k)\to v\not\in {\rm dom} \nabla f$.

From $\nabla f(a_{k,\lambda}) = \nabla f(b_k) + \lambda (\nabla f(a_k)-\nabla f(b_k))$ and boundedness of the $\nabla f(a_k)$
follows $\|\nabla f(a_{k,\lambda})\| \propto \lambda$ uniformly over $k$. 
Since $\|\nabla f(a_{k,\lambda})\|\to \infty$ for fixed $k$, $a_{k,\lambda}$ must be defined for all $\lambda > 0$. 
Now parametrize the same perpendicular curve as $\nabla f(\bar{a}_{k,\mu})= \nabla f(b_k) + \mu (\nabla f(\bar{a}_k)-\nabla f(b_k))$, then 
again $\nabla f(\bar{a}_{k,\mu}) \propto \mu$ uniformly over $k$. With the same argument as above, $\bar{a}_{k,\mu}$ must be defined for all
$\mu > 0$.

Now for every $k$ and $\lambda > 0$ there exists $\mu_k(\lambda)$ such that $a_{k,\lambda} = \bar{a}_{k,\mu_k(\lambda)}$. 
But the relation between the two is given by
\begin{equation}
\label{mix}
\lambda (\nabla f(a_k)-\nabla f(b_k)) = \mu( \nabla f(\bar{a}_k)-\nabla f(b_k)),
\end{equation}
hence by boundedness of the $\nabla f(a_k)$ we have  $r\lambda \leq \mu_k(\lambda) \leq r'\lambda$ for all $k$ with certain $r,r'$
not depending on $k$.  

Since $\bar{a}_k$ is between $a_k$ and $b_k$, we have $\mu > 1$ when $\lambda=1$. We argue that this implies $\mu_k(\lambda) > \lambda$ for almost all $k$. Indeed,
if for some $k$ there is a moment, where $\mu_k(\lambda) = \lambda$, then $\nabla f(a_k) = \nabla f(\bar{a}_k)$, and that could happen only a finite number of times, because
$\nabla f(\bar{a}_k) \to \nabla f(a^*)$, while $\nabla f(a_k)\to v\not\in {\rm dom} \nabla f$. Hence we can assume $\mu_k(\lambda) > \lambda$ for all $k$. 
From (\ref{mix}) we now get
$$
\left(1-\frac{\lambda}{\mu_k(\lambda)}\right) \nabla f(b_k) + \frac{\lambda}{\mu_k(\lambda)} \nabla f(a_k) = \nabla f(\bar{a}_k),
$$
a convex combination, 
and passing to the limit $(k\to \infty)$ in a subsequence $\mu_k(\lambda) \to \mu^* \geq \lambda$,
$$
\left(1-\frac{\lambda}{\mu^*}\right) \nabla f(b^*) + \frac{\lambda}{\mu^*} v = \nabla f(a^*),
$$
proving $\nabla f(a^*)\in [\nabla f(b^*),v]$.
But  the $\bar{a}_k$
are independent of $v$, so we can arrange the same estimate for other points ${a}^*$  on the limit perpendicular $a^*_\lambda$ generated by $b^*\in \cev{P}_B(a^*)$.
But that means $\nabla f(a^*_\lambda) \in [\nabla f(b^*),v]$ for the entire perpendicular, contradicting $\|\nabla f(a^*_\lambda)\| \to \infty$. 
\end{proof}

\section{Gaps between sets}
\label{sect_gap}
Let $a_k,b_k$ be a Bregman alternating sequence. Since
$D(b_k,a_k) \leq D(b_{k-1},a_k) \leq D(b_{k-1},a_{k-1})$, monotone convergence gives $D(b_k,a_k) \to \frac{1}{2} r^{*2}$, and
$D(b_{k-1},a_k)\to \frac{1}{2}r^{*2}$ for some $r^*\geq 0$.
Now let $k\in N$ be an infinite subsequence of $\mathbb N$ such that $b_{k-1} \to b^*$, $a_k\to a^*$, $b_k \to \hat{b}$. 
As $a_k\in \vec{P}_A(b_{k-1})$, we have $a^*\in \vec{P}_A(b^*)$. Similarly,
as
$b_k \in \cev{P}_B(a_k)$, we have $\hat{b}\in \cev{P}_B(a^*)$. But $D(b^*,a^*) = D(\hat{b},a^*) = \frac{1}{2}r^{*2}$, hence
$b^*\in \cev{P}_B(a^*)$, too. So we have found a pair $(b^*,a^*) \in B \times A$ with $a^*\in \vec{P}_A(b^*)$,
$b^*\in \cev{P}_B(a^*)$. We write $b^* \sim a^*$ for such pairs. 

Let $A^*,B^*$ be the sets of accumulation points of the $a_k,b_k$. The above argument
shows that for every $a^*\in A^*$ there exists $b^*\in B^*$ such that $b^* \sim a^*$, and
 for every $b^*\in B^*$ there exists $a^*\in A^*$ with $b^* \sim a^*$. We call
$K^* = \{(b^*,a^*)\in B^* \times A^*: b^* \sim a^*\}$  {\it the gap} of the alternating sequence.  The case $r^*=0$ is not excluded,
where of course $b^*\sim a^*$ implies $b^*=a^*\in A \cap B$.

Abstracting from the sequence $a_k,b_k$, we call a pair
$(K^*,r^*)$ a {\it gap} between $A$ and $B$ if $K^*\subset B\times A$ is compact
with $b^*\sim a^*$ for all $(b^*,a^*)\in K^*$ and $D(b^*,a^*)=\frac{1}{2}r^{*2}$.
Note that every $x^*\in A \cap B$ gives rise to a zero gap $(\{(x^*,x^*)\},0)$.

\begin{lemma}
\label{lem5}
Let $F(x,y) = i_B(x) + D(x,y) + i_A(y)$, and suppose
$(b^*,a^*)\in K^*$, then
$(0,0) \in \partial F(b^*,a^*)$.
\end{lemma}

\begin{proof}
Consider a building block $b \stackrel{r}{\longrightarrow} a^+ \stackrel{l}{\longrightarrow} b^+$, then
$n_B(b^+) = \nabla f(a^+)-\nabla f(b^+)\in N_B(b^+)$ and $n_A(a^+) = \nabla^2f(a^+)(b-a^+)\in N_A(a^+)$ by Lemma \ref{prox_normals}, hence
\begin{equation}
\label{oefter}
(\lambda n_B(b^+) + \nabla f(b^+)-\nabla f(a^+),\mu n_A(a^+) + \nabla^2f(a^+)(a^+-b^+)) \in \partial F(b^+,a^+)
\end{equation}
for all $\lambda,\mu \geq 0$.
Since $b^*\sim a^*$,  we may choose $b=b^+=b^*$ and $a^+=a^*$ in the building block, which gives $(0,0) \in \partial F(b^*,a^*)$ on putting
$\lambda = \mu = 1$.
\end{proof}

\begin{proposition}
\label{prop_compact}
Suppose $A,B,f$ are definable. Then there is only a finite number of gap values
$0 \leq r_1 < r_2 < \dots < r_N$. There exist $\eta_i > 0$, $i=1,\dots,N$, such that
every alternating sequence $a_k,b_k$ which satisfies $\frac{1}{2}r_i^2\leq D(b_k,a_k) < \frac{1}{2}r_i^2 + \eta_i$ must have value convergence
$D(b_k,a_k)\to \frac{1}{2}r_i^2$, $D(b_{k-1},a_k) \to \frac{1}{2}r_i^2$.
\end{proposition}

\begin{proof}
This follows with \cite[Prop. 2]{kurdyka}.
\end{proof}

\section{Angle condition}
\label{sect_angle_rl}
In this section we  introduce the angle condition and show that it
is a geometric form of the Kurdyka-\L ojasiewicz inequality.

\begin{definition}
\label{def_angle_gap}
{\bf 
(Angle condition)}.
{\rm 
Let $\sigma:(0,\infty) \to (0,\infty)$ be increasing. 
The set $B$ satisfies the $rl$-angle condition with constant $\gamma$ and shrinking function $\sigma$ with respect to $A$ at a gap pair $b^* \sim a^*$ with gap value $r^*$, if
there exists a neighborhood $W$ of $(b^*,a^*)$ and $\eta > 0$ such that
\begin{equation}
\label{angle}
\frac{1-\cos \alpha}{\sigma(\cev{D}_B(a^+)-\textstyle\frac{1}{2}r^{*2})} \geq \gamma
\end{equation}
for every building block $b \stackrel{r}{\longrightarrow} a^+ \stackrel{l}{\longrightarrow} b^+$ with $(b^+,a^+)\in W$, 
$\frac{1}{2}r^{*2} \leq D(b^+,a^+) < \frac{1}{2}r^{*2} + \eta$, where
$\alpha = \angle(b-a^+,b^+-a^+)$.
}
\end{definition}

\begin{remark}
A standard compactness argument shows that if $(K^*,r^*)$ is a gap
such that the angle condition holds at every $(b^*,a^*)\in K^*$, then it holds for all building blocks in a neighborhood of $K^*$
with the same $\sigma$, $\eta,\gamma$.
\end{remark}

We now show that the angle condition can be understood as  a geometric form of the K\L-inequality.

\begin{proposition}
\label{prop_gap_angle}
Let $(K^*,r^*)$ be a gap between $A$ and $B$.
Suppose $F(x,y) = i_B(x) + D(x,y)+ i_A(y)$ satisfies the K\L-condition at every $(b^*,a^*)\in K^*$.
Then the angle condition is satisfied in the following form. There exists a neighborhood $W$
of $K^*$, a de-singularizing function $\phi$, and constants $\gamma,\eta > 0$ such that every building block
$b \stackrel{r}{\longrightarrow} a^+ \stackrel{l}{\longrightarrow} b^+$ with $(b^+,a^+)\in W$ 
and $\frac{1}{2}r^{*2} \leq D(b^+,a^+) < \frac{1}{2}r^{*2} + \eta$ satisfies
\begin{equation}
\label{angle_from_KL}
\phi'(\cev{D}_B(a^+)-\textstyle\frac{1}{2}r^{*2})^2 \cev{D}_B(a^+)(1-\cos\alpha) \geq \gamma,
\end{equation}
where 
$\alpha = \angle (b-a^+,b^+-a^+)$.
\end{proposition}

\begin{proof}
$K^*$ being compact, and $F$ having constant value $\frac{1}{2}r^{*2}$
on $K^*$,  the K\L-inequality is satisfied as follows: There exists a bounded neighborhood $W$ of $K^*$, a de-singularizing function $\phi$, and constants $\gamma,\eta > 0$
such that 
$$
\phi'(F(b^+,a^+) - \textstyle\frac{1}{2}r^{*2}) {\rm dist}_{|\cdot|}((0,0),\partial F(b^+,a^+)) \geq \gamma
$$
for all $(b^+,a^+) \in W$ with $\frac{1}{2}r^{*2} \leq F(b^+,a^+) < \frac{1}{2}r^{*2} + \eta$, and where $|(x,y)| = \|x\|+\|y\|$.

Now consider building blocks  $b \stackrel{r}{\longrightarrow} a^+ \stackrel{l}{\longrightarrow} b^+$ with $(b^+,a^+)\in W$ and $\frac{1}{2}r^{*2} \leq D(b^+,a^+) < \frac{1}{2}r^{*2} + \eta$.  
Since
$$
(\lambda n_B(b^+) + \nabla f(b^+)-\nabla f(a^+),\mu n_A(a^+)+\nabla^2f(a^+)(a^+-b^+))\in \partial F(b^+,a^+)
$$
for all $n_B(b^+) \in N_B(b^+)$ and $n_A(a^+)\in N_A(a^+)$, 
Lemma \ref{prox_normals} gives
\begin{align*}
\phi'(F(b^+,a^+) - \textstyle\frac{1}{2}r^{*2}) \left(
\|\lambda (\nabla f(a^+)-\nabla f(b^+)) + \nabla f(b^+)-\nabla f(a^+))\| \right.
+ \qquad \\
\left.
\|\mu \nabla^2f(a^+)(b-a^+) + \nabla^2f(a^+)(a^+-b^+)\| \right)\geq \gamma
\end{align*}
for all $\lambda,\mu \geq 0$. Choosing $\lambda = 1$, 
gives
$$
\phi'(F(b^+,a^+) - \textstyle\frac{1}{2}r^{*2})  \|\mu \nabla^2f(a^+)(b-a^+) + \nabla^2f(a^+)(a^+-b^+)\| \geq \gamma 
$$
for all $\mu \geq 0$. Now due to boundedness of $W$  the eigenvalues of all Hessians $\nabla^2f(a^+)$
with $(b^+,a^+)\in W$ are bounded above by a constants $\Lambda$.
Therefore  $\|\mu \nabla^2f(a^+)(b-a^+) + \nabla^2f(a^+)(a^+-b^+)\| \leq \Lambda  \|\mu (b-a^+) + a^+-b^+\|$, so that
$$
\phi'(F(b^+,a^+) - \textstyle\frac{1}{2}r^{*2})  \|\mu (b-a^+) + a^+-b^+\| \geq \Lambda^{-1}\gamma 
$$
for all $\mu \geq 0$. Passing to the infimum over $\mu \geq 0$ implies
$$
\phi'(F(b^+,a^+) - \textstyle\frac{1}{2}r^{*2}) \sin\alpha \|a^+-b^+\| \geq \Lambda^{-1}\gamma,
$$
for angles $\alpha = \angle (b-a^+,b^+-a^+)$ less than $90^\circ$, while for angles larger than $90^\circ$ the minimum is attained at 
$\phi'(F(b^+,a^+) - \textstyle\frac{1}{2}r^{*2})  \|a^+-b^+\|$.  As the statement is clear in the latter case, we continue with $\alpha < 90^\circ$. 
Here  using  $\|a^+-b^+\| \leq M D(b^+,a^+)^{1/2}$, taking squares gives
$$
\phi'(\cev{D}_B(a^+)-\textstyle\frac{1}{2}r^{*2})^2 \cev{D}_B(a^+) \sin^2\alpha \geq \gamma^2/M^2 \Lambda^2
$$
and then (using $1-\cos \alpha \geq \frac{1}{2}\sin^2\alpha$):
$$
\phi'(\cev{D}_B(a^+)-\textstyle\frac{1}{2}r^{*2})^2 \cev{D}_B(a^+) (1-\cos\alpha) \geq \gamma^2/2M^2\Lambda^2 =:\gamma',
$$
where $\gamma'$ depends only on $K^*$ and the bounds (\ref{bounds})  associated with it.
This proves the claim.
\end{proof}

\begin{remark}
Clearly (\ref{angle_from_KL}) gives the $rl$-angle condition (\ref{angle}) with
angle shrinking function  $\sigma(s) = \phi'(s)^{-2} s^{-1}$ when $r^*=0$. When $r^* > 0$, then the term $\cev{D}_B(a^+)$  stays bounded away from 0
and hence does not contribute to the singularity. We can then shuffle it to the right, obtaining yet another $\gamma''$ depending only on $K^*$
and $r^*$, now with $\sigma(s) = \phi'(s)^{-2}$. During the following we will always use the angle condition with
angle shrinking functions generated by de-singularizing functions $\phi$. 
\end{remark}

\begin{corollary}
Let $A,B,f$ be definable. Then there exists a unique de-singularizing function $\phi$
which works for each of the finitely many gap values. The angle shrinking function
is $\sigma(s) = \phi'(s)^{-2}$ for gaps $r_i > 0$, and $\sigma(s) = \phi'(s)^{-2}s^{-1}$ for $r_1$  in case $r_1 = 0$. 
\end{corollary}

\subsection{Tangential and transversal intersection}
\label{sect_tangential}

 It is clear that
the mutual geometric position of $A,B$ near a point $\bar{x}\in A \cap B$ must be decisive for the speed of convergence
of alternating Bregman projections near $\bar{x}$. We distinguish
tangential and 
transversal intersection, expecting transversality  to give linear speed, while tangential intersection should force a slowdown to sub-linear speed. 
We now give an interpretation of these using the angle condition.

\begin{definition}
{\bf (Transversality)}.
{\rm
We say that $B$ intersects $A$ {\it rl-transversally} at $\bar{x}\in A \cap B$, 
  if there exists a neighborhood $U$
of $\bar{x}$ and $\underline{\alpha} > 0$ such that 
for every building block $b \stackrel{r}{\longrightarrow} a^+ \stackrel{l}{\longrightarrow} b^+$ in $U$
the angle 
$\alpha = \angle(b-a^+,b^+-a^+)$
is larger than $\underline{\alpha}$. Otherwise we say that
$B$ intersects $A$ {\it rl-tangentially} at $\bar{x}$. 
}
\end{definition}

This notion is weaker than classical transversality, and yet guarantees linear convergence
of the alternating method, as will indeed be proved in Corollary \ref{cor_transversal}. One classical notion of transversality in non-smooth calculus is the following:

\begin{proposition}
\label{transversal}
Suppose $N_B(\bar{x}) \cap (-N_A(\bar{x}) )= \{0\}$. Then $B$ intersects $A$ rl-transversally at $\bar{x}\in A \cap B$.
\end{proposition}

\begin{proof}
Assume on the contrary that there exist
building blocks
$b_{k-1} \stackrel{r}{\longrightarrow} a_k \stackrel{l}{\longrightarrow} b_k$ 
with $b_{k-1},a_k,b_k \to \bar{x}$ as $k \to \infty$ such that 
$\alpha_k = \angle(b_{k-1}-a_k,b_k-a_k)\to 0$. Let $u_k = (b_{k-1}-a_k)/\|b_{k-1}-a_k\|$, $v_k = (b_k-a_k)/\|b_k-a_k\|$, 
and
$w_k=(\nabla f(a_k)-\nabla f(b_k) )/ \|a_k-b_k\|$.

Recall 
$n_A(a_k)=\nabla^2f(a_k)(b_{k-1}-a_k) \in N_A(a_k)$ and  $n_B(b_k) = \nabla f(a_k)-\nabla f(b_k) \in N_B(b_k)$ by Lemma \ref{prox_normals}.
Hence $\nabla^2f(a_k)u_k \in N_A(a_k)$ and $w_k\in N_B(b_k)$.
Select an infinite subsequence $k\in K$ of $\mathbb N$ such that $u_k\to u \not=0$, $v_k\to v\not=0$, $w_k\to w\not=0$, 
the latter due to (\ref{nabla_bounds}). 
Then
$\nabla^2f(a_k)u_k \to \nabla^2f(\bar{x})u \in N_A(\bar{x})$, and $\nabla^2 f(a_k)v_k \to \nabla^2f(\bar{x})v$,
$w_k \to w\in N_B(\bar{x})$. But $\angle (u_k,v_k) \to 0$ implies $\angle(u,v)=0$, and 
since $u,v$ are unit vectors, we have $u=v$.  Since $\nabla^2f(\bar{x})u\in N_A(\bar{x})$, this gives $\nabla^2f(\bar{x})v \in N_A(\bar{x})$. 

Taylor-Young expansion of $\nabla f$ at the $a_k$
gives
$\nabla f(a_k) - \nabla f(b_k) = \nabla^2f(a_k)(a_k-b_k) + o(\|a_k-b_k\|)$, where the o-term may be made uniform on a bounded neighborhood of $\bar{x}$, see Section \ref{sect_uTJ}. 
Therefore 
$$
w_k=\frac{\nabla f(a_k)-\nabla f(b_k)}{\|a_k-b_k\|} = \nabla^2 f(a_k) (-v_k) + o(1) \to \nabla^2 f(\bar{x})(-v) \in -N_A(\bar{x}).
$$
Since $w_k \to w\in N_B(\bar{x})$, we have $w\in N_B(\bar{x}) \cap (-N_A(\bar{x}))$, a contradiction.
\end{proof}

\begin{remark}
In \cite{bauschke1,bauschke2} the authors propose a refinement in the euclidean case, where $N_B(\bar{x})$ is replaced by a restricted
normal cone $N_B^A(\bar{x})$, which only considers projections on $B$ stemming from $A$, and similarly for $N_A^B(\bar{x})$. This could be extended 
to left and right Bregman projections. We do not pursue this natural idea further. The proof in
\cite[Prop. 1]{noll} may be adapted to the Bregman setting.
\end{remark}

When $B$ intersects $A$ $rl$-tangentially at $\bar{x}$, then there are $rl$-building blocks  for which the angle $\alpha=\angle(b-a^+,b^+-a^+)$ shrinks to $0$ as these
$b \stackrel{r}{\longrightarrow} a^+ \stackrel{l}{\longrightarrow} b^+$ approach
$\bar{x}$. 
Following \cite{noll,gerchberg}, we now relate this to the angle condition
in the zero gap case $r^*=0$. Note that for $r^*=0$ condition (\ref{angle}) means
\begin{equation}
\label{angle_no}
\frac{1-\cos \alpha}{\sigma(\cev{D}_B(a^+))} \geq \gamma
\end{equation}
for every building block $b \stackrel{r}{\longrightarrow} a^+ \stackrel{l}{\longrightarrow} b^+$ in $U$, where
$\alpha = \angle(b-a^+,b^+-a^+)$. Therefore when the angle $\alpha$ tends to $0$, meaning $1-\cos\alpha \to 0$, then
the speed with which it is allowed to do so is controlled by the angle shrinking function
$\sigma$, which when obtained from the K\L-inequality is $\sigma(s) =  \phi'(s)^{-2}s^{-1}$ for a de-singularizing function $\phi$.
This was introduced in \cite{noll,gerchberg} for the euclidean case with shrinking functions
$\sigma(s) = s^\theta$ for $\theta\in [\frac{1}{2},1)$.

When (\ref{angle_no}) holds but $\sigma(\cev{D}_B(a^+))$ does {\it not} shrink to $0$, then $\alpha$ must also stay away from $0$, which is precisely when 
$B$ intersects $A$  $rl$-transversally at $\bar{x}$. Hence this case is covered by the angle condition.

\begin{definition}
{\bf (Dual transversality)}.
{\rm 
We say that $A$ intersects $B$ {\it lr-transversally} at $\bar{x}\in A \cap B$ if there exists a neighborhood $U$ of $\bar{x}$ and $\underline{\alpha} > 0$ such that
for every building block $a\stackrel{l}{\longrightarrow} b \stackrel{r}{\longrightarrow} a^+$ in $U$ the angle $\alpha = \angle (\nabla f(a)-\nabla f(b),\nabla f(a^+)-\nabla f(b))$ is larger than $\underline{\alpha}$.
Otherwise $A$ intersects $B$ $lr$-tangentially at $\bar{x}$.
}
\end{definition}

As a first application of duality we have

\begin{proposition}
Suppose $N_B(\bar{x}) \cap (-N_A(\bar{x})) = \{0\}$. Then $A$ intersects $B$ $lr$-transversally at $\bar{x}\in A \cap B$.
\end{proposition}

\begin{proof}
Let $A^*=\nabla f(B)$, $B^*=\nabla f(A)$, then 
$\nabla^2f(\bar{x}) {N}_{B^*}(\nabla f(\bar{x}) )= {N}_A(\bar{x})$ by the chain rule \cite[Thm. 10.6, Ex. 6.7]{rock}, applied to the $C^1$-diffeomorphism $\nabla f$.
Similarly $\nabla^2f(\bar{x}){N}_{A^*}(\nabla f(\bar{x})) ={N}_B(\bar{x})$. Since $\nabla^2f(\bar{x}) \succ 0$,
the hypothesis implies
$N_{A^*} (\nabla f(\bar{x})) \cap (-N_{B^*} (\nabla f(\bar{x}))) = \{0\}$. 
Therefore by the previous proposition $B^*$ intersects $A^*$ $rl$-transversally at $\nabla f(\bar{x}) \in A^* \cap B^*$.
For 
building blocks 
$b^* \stackrel{r*}{\longrightarrow} a^{*+} \stackrel{l}{\longrightarrow} b^{*+}$ close to $\nabla f(\bar{x})$ this means
$\angle (b^*-a^{*+},b^{*+}-a^{*+}) \geq \underline{\alpha} > 0$. 
But
$rl$-building blocks
$a\stackrel{l}{\longrightarrow} b \stackrel{r}{\longrightarrow} a^+$
in the neighborhood of $\bar{x}$
are mapped by $\nabla f$ to 
$lr$-building block 
$b^* \stackrel{r*}{\longrightarrow} a^{*+} \stackrel{l}{\longrightarrow} b^{*+}$
in the neighborhood of $\nabla f(\bar{x})$
and this gives directly $\alpha = \angle (\nabla f(a)-\nabla f(b),\nabla f(a^+)-\nabla f(b)) \geq \underline{\alpha}$.
\end{proof}

This calls for a dual angle condition, which we will obtain shortly.

\subsection{Duality for the K\L-inequality}
\label{sect_dual_kl}
We show that the K\L-inequality is directly amenable to duality. 
\begin{lemma}
\label{dual_kl}
Let $F(x,y) = i_B(x) + D(x,y) + i_A(y)$ and put
$F_*(u,v) = i_{B^*}(u) + D^*(u,v) + i_{A^*}(v)$, where $A^*=\nabla f(B)$, $B^*=\nabla f(A)$. Suppose
$F$ satisfies a K\L-inequality at $(\bar{x},\bar{y})$, then $F_*$
satisfies a K\L-inequality at $(\bar{u},\bar{v}) = (\nabla f(\bar{y}),\nabla f(\bar{x}))$ with the same de-singularizing function.
\end{lemma}

\begin{proof}
We have $F(x,y) 
= i_{B\times A}(x,y) + D(x,y)$, hence $\partial F(x,y) = N_{B\times A}(x,y) + \partial D(x,y)
= N_B(x) \times N_A(y) + \partial D(x,y)$, $D$ being jointly differentiable, and using \cite[Prop. 6.41]{rock}. Now $D(x,y) = D^*(\nabla f(y),\nabla f(x))$, and 
$i_A(y) = i_{\nabla f(A)}(\nabla f(y)) = i_{B^*}(\nabla f(y))$,
$i_B(x) = i_{A^*}(\nabla f(x))$. Hence 
$$
F(x,y) = F_*(\Phi(x,y))
$$
under the isomorphism $\Phi(x,y) = (\nabla f(y),\nabla f(x))$,
where $F_*(u,v) = i_{B^*\times A^*}(u,v) + D^*(u,v)$. Then $\partial F(x,y) = \Phi'(x,y)^{T} \partial F_*(\Phi(x,y))$ by the chain rule \cite[Thm. 10.6]{rock}.  
Suppose $(0,0)\in \partial F(\bar{x},\bar{y})$, then as $F$ satisfies the K\L-condition at $(\bar{x},\bar{y})$, we have
$$
\phi'(F(x,y)-F(\bar{x},\bar{y})) {\rm dist}((0,0),\partial F(x,y)) \geq \gamma
$$
for some $\gamma > 0$, some $\eta > 0$, a neighborhood $W$ of $(\bar{x},\bar{y})$, and every $(x,y)\in W$
satisfying $F(\bar{x},\bar{y}) \leq F(x,y) \leq F(\bar{x},\bar{y}) + \eta$. Now  let $u = \nabla f(y)$, $v = \nabla f(x)$,
$\bar{u} = \nabla f(\bar{y})$, $\bar{v}=\nabla f(\bar{x})$. Then
$(0,0) \in \partial F_*(\nabla f(\bar{y}),\nabla f(\bar{x}))$, because $\Phi'(\bar{x},\bar{y})^{T}$ is a linear isomorphism.
Let $V = \Phi(W)$, then $V$ is a neighborhood of $(\bar{u},\bar{v})$. Let $(u,v) \in V$, $(u,v) = \Phi(x,y)$ with $(x,y)\in W$.
Let $h \in \partial F_*(u,v)$, then $h = \Phi'(x,y)^{-T} g$ for $g\in \partial F(x,y)$. But then
\begin{align*}
\phi'(F_*(u,v)-F_*(\bar{u},\bar{v})) \|h\| &= \phi'(F(x,y)-F(\bar{x},\bar{y})) \| \Phi'(x,y)^{-T} g\|\\& \geq \phi'(F(x,y)-F(\bar{x},\bar{y})) k \|g\| \geq k\gamma,
\end{align*}
where we use $\|\Phi'(x,y)^{T}\| \leq k^{-1}$ for all $(x,y)\in W$. Here $\Phi'(x,y) = {\rm diag}(\nabla^2f(y),\nabla^2f(x)) = \Phi'(x,y)^T$,
and we can without loss choose $W = W_1 \times W_2$ such that $\nabla^2f(y)$ is bounded away from $0$ on the neighborhood $W_2$
of $\bar{y}$, $\nabla^2f(x)$ bounded away from $0$ on the neighborhood $W_1$ of $\bar{x}$.
This means, the K\L-inequality for $F_*$ is satisfied with $\gamma'=k\gamma$ and the same $\phi$ and $\eta$.
\end{proof}

\subsection{Duality for the angle condition}
\label{sect_5.3}
We now show that duality leads the way to the correct dual angle condition. Observe that
$\sigma(\cev{D}_B(a^+)-\frac{1}{2}r^{*2}) = \sigma(\vec{D}_{A^*}(b^*)-\frac{1}{2}r^{*2})$ due to (\ref{fenchel}) and (\ref{dual_formula}).
Also $\alpha = \angle (b-a^+,b^+-a^+)$ is the same as
$\alpha = \angle (\nabla f^*(a^*)-\nabla f^*(b^{*+}),\nabla f^*(a^{*+})-\nabla f^*(b^{*+}))$. This means the correct definition is as follows:
\begin{definition} 
\label{def_angle_II}
{\bf (Dual angle condition)}.
{\rm
Let $\sigma:(0,\infty) \to (0,\infty)$ be monotonically increasing. The set $A$ satisfies the $lr$-angle condition with constant $\gamma > 0$
and shrinking function $\sigma$
with respect to $B$ at a gap pair $b^* \sim a^*$, if there exists a neighborhood $W$ of $(b^*,a^*)$ and $\gamma,\eta > 0$ such that
\begin{equation}
\label{angle_lr}
\frac{1-\cos \alpha}{\sigma(\vec{D}_A(b)-\frac{1}{2}r^{*2})} \geq \gamma
\end{equation}
for every $lr$-building block  $a\stackrel{l}{\longrightarrow} b \stackrel{r}{\longrightarrow} a^+$ with $(b,a^+)\in W$ and
$\frac{1}{2}r^{*2} \leq \vec{D}_A(b) < \frac{1}{2}r^{*2} + \eta$, where
$\alpha = \angle(\nabla f(a)-\nabla f(b),\nabla f(a^+)-\nabla f(b))$.}
\end{definition}

This can also be extended to gaps $(K^*,r^*)$ using duality. Again we will content ourselves with angle shrinking functions of the form
$\sigma(s) = \phi'(s)^{-2}s^{-1}$, respectively, $\sigma(s) = \phi'(s)^{-2}$ for a de-singularizing $\phi$.

\subsection{Local projections}
\label{sect_local}
The angle condition has an advantage over the K\L-condition. 
We say that $a_k\in \vec{P}_A(b_k)$ is a local projection
when the minimum (\ref{right}) is local, 
while still guaranteeing decrease $D(b_k,a_{k-1}) > D(b_k,a_k)$. 
Consider a bounded alternating sequence which is local in this sense, in symbols $b_{k-1} \stackrel{\,r\,}{\dashrightarrow} a_k \stackrel{l}{\longrightarrow} b_k$,   
and let $A^*,B^*$ be the set of accumulation points of the $a_k,b_k$.  Define
$A^s = \{a_k: k\in \mathbb N\} \cup A^*$, $B^s = \{b_k: k\in \mathbb N\}\cup B^*$. Then $a_k,b_k$ is a standard
alternating sequence $b_{k-1} \stackrel{r}{\longrightarrow} a_k \stackrel{l}{\longrightarrow} b_k$ between the closed sets $A^s,B^s$, as points in $A$ which previously might have made $\vec{P}_A$ local have been removed. 

Now suppose $A,B,f$ are definable, then the K\L-inequality
holds everywhere. Here we use the fact that the proof of Proposition \ref{prop_gap_angle} is still valid, since
$\partial F(b^+,a^+)$ is not altered by $\vec{P}_A$ being local. 
Hence the angle condition holds everywhere. But the angle condition is expressed in terms of building blocks,
hence it is untouched when we restrict to the smaller sets $A^s,B^s$. This is significant, because
$A^s,B^s$ have no reason to be definable. Nor is there a good argument in favor of
$F^s(x,y) = i_{B^s}(x) + D(x,y) + i_{A^s}(y)$  still having the K\L-property, because the normal cones $N_{A^s}, N_{B^s}$ are larger than $N_A,N_B$. 
This thread will be picked up again in Corollary \ref{cor_local}.

\section{Three point inequality and duality}
\label{sect_three_point}
In this section we discuss the second fundamental ingredient of our convergence theory, referred to as the three-point inequality in \cite{noll}.
In the convex setting, a related notion has been discussed in \cite{csiszar1}. We extend this now to Bregman projections.

\begin{definition}
{\rm 
The building block
$b \stackrel{r}{\longrightarrow} a^+ \stackrel{l}{\longrightarrow} b^+$ satisfies the $rl$-three-point inequality with
constant $\ell \in (0,1]$ if
\begin{equation}
\label{rlr_3_pt}
D(b,a^+) \geq D(b^+,a^+) + \ell D(b,b^+).
\end{equation}
The building block
$a \stackrel{l}{\longrightarrow} b \stackrel{r}{\longrightarrow} a^+$ satisfies the $lr$-three-point inequality with $\ell \in (0,1]$
if
\begin{equation}
\label{lrl_3_pt}
D(b,a) \geq D(b,a^+) + \ell D(a^+,a).
\end{equation}
}
\end{definition}

\begin{remark}
Putting $b^*=\nabla f(a^+)$, $a^* = \nabla f(b)$, $a^{*+} = \nabla f(b^+)$, 
we see from the dual relationship $D(x,y) = D^*(\nabla f(y),\nabla f(x))$ that (\ref{rlr_3_pt}) transforms into
$$
D^*(b^*,a^*) \geq D^*(b^*,a^{*+}) + \ell D^*(a^{*+},a^*),
$$ 
which is (\ref{lrl_3_pt}) for building blocks $a^* \stackrel{l*}{\longrightarrow} b^* \stackrel{r*}{\longrightarrow} a^{*+}$
alternating between
$B^* = \nabla f(A)$ and $A^* = \nabla f(B)$ due to (\ref{dual_formula}). In other words, the three-point inequality is directly amenable to
duality.
\end{remark}
In the euclidean case \cite{noll} 
$lr$- and $rl$-variants coincide.  Now we define:

\begin{definition}
{\rm 
Let $a^*\in A$, $b^*\in B$, $b^* \sim a^*$. We say that the $rl$-three-point inequality holds at  $(b^*,a^*)$ 
if there exists $\ell \in (0,1)$ and $\delta > 0$ such that every $rl$-building block $b \stackrel{r}{\longrightarrow} a^+ \stackrel{l}{\longrightarrow} b^+$
with $b^+ \in B(b^*,\delta)$, $a^+\in B(a^*,\delta)$  
satisfies the $rl$-three-point inequality with constant $\ell$.
}
\end{definition}

The definition for the $lr$-case is analogous.

\begin{proposition}
\label{prop_convex}
Suppose $B$ is convex. Then the $rl$-three-point inequality holds for all $A$. Suppose 
$\nabla f(A)$ is convex. Then the $lr$-three-point inequality holds for all $B$.   
\end{proposition}

\begin{proof}
By the Bregman law of cosines $D(b,a^+)=D(b^+,a^+) + D(b,b^+) - \langle \nabla f(a^+)-\nabla f(b^+),b-b^+\rangle$.
Now $\nabla f(a^+)-\nabla f(b^+) \in N_B(b^+)$, hence  by convexity of $B$ the term $- \langle \nabla f(a^+)-\nabla f(b^+),b-b^+\rangle$ is positive for $b\in B$,
and we get $D(b,a^+)\geq D(b^+,a^+) + D(b,b^+)$, which is (\ref{rlr_3_pt}) with $\ell=1$. For the second statement, by duality
we have to show that building blocks $b^*\stackrel{r*}{\longrightarrow}a^{*+}\stackrel{l*}{\longrightarrow}b^{*+}$
satisfy the  $rl$-three-point inequality. This follows from the first part, as now $B^*=\nabla f(A)$ is convex, while $A^*=\nabla f(B)$.
\end{proof}

\begin{proposition}
\label{prop7}
Let $(K^*,r^*)$ be a gap between $A$ and $B$, and suppose
that
at every $(b^*,a^*)\in K^*$ a three-point inequality for rl-building blocks is satisfied. Then there exists $0 < \ell < 1$ and a neighborhood $W$ of $K^*$ such that the three point inequality holds with
the same $\ell\in (0,1]$  for all building blocks $b \stackrel{r}{\longrightarrow} a^+ \stackrel{l}{\longrightarrow} b^+$
satisfying  $(b^+,a^+) \in W$. 
\end{proposition}

\begin{proof}
This is a compactness argument. By hypothesis every $(b^*,a^*)\in K^*$ has a neighborhood $B(b^*,\delta_{b^*,a^*}) \times B(a^*,\delta_{b^*,a^*})$
such that $rl$-building blocks with $b^+\in B(b^*,\delta_{b^*,a^*})$,
$a^+\in B(a^*,\delta_{b^*,a^*})$ satisfy the three-point inequality with $\ell_{b^*,a^*}\in (0,1)$. Now 
$K^*\subset \bigcup\{ B(b^*,\delta_{b^*,a^*}/2) \times B(a^*,\delta_{b^*,a^*}/2): (b^*,a^*)\in K^*\}$, and by compactness of $K^*$
there exist finitely many pairs $(b^*_i,a^*_i)\in K^*$ such that
$K^*\subset \bigcup\{ B(b^*,\delta_{b^*_i,a^*_i}/2) \times B(a^*,\delta_{b^*_i,a^*_i}/2): i=1,\dots,m\} =:W$.
Now let $\ell = \min\{\ell_{b_i^*,a_i^*}: i=1,\dots,m\}$, 
and let $b \stackrel{r}{\longrightarrow} a^+ \stackrel{l}{\longrightarrow} b^+$
be a building block with $(b^+,a^+)\in W$. 
Then $(b^+,a^+) \in  B(b^*,\delta_{b^*_i,a^*_i}/2) \times B(a^*,\delta_{b^*_i,a^*_i}/2)$
for some $i$. 
Therefore this building block satisfies the three-point inequality with $\ell_{b_i^*,a_i^*}$, hence also with $\ell$.
\end{proof}

Interestingly, the three point inequality for building blocks $b \stackrel{r}{\longrightarrow} a^+ \stackrel{l}{\longrightarrow} b^+$ approaching a gap forces $b,b^+$
to get close.

\begin{lemma}
\label{lem7}
Let $a_k,b_k$ be a Bregman alternating sequence with gap $(K^*,r^*)$. Suppose the rl-three point inequality holds in a neighborhood of $K^*$. 
Then $b_{k-1}-b_k \to 0$.
\end{lemma}

\begin{proof}
Select an infinite subsequence $k\in N$ such that $b_{k-1} \to b^*$, $a_k\to a^*$, $b_k \to \hat{b}$. Suppose $\hat{b} \not=b^*$. We have
$D(b_{k-1},a_k) \to \frac{1}{2}r^{*2}$, $D(b_k,a_k) \to \frac{1}{2} r^{*2}$, $D(b_{k-1},b_k) \to D(b^*,\hat{b}) > 0$. Therefore
$D(b_{k-1},a_k) \geq  D(b_k,a_k) + \ell D(b_{k-1},b_k) \to \frac{1}{2}r^{*2} + \ell D(b^*,\hat{b}) > \frac{1}{2}r^{*2}$, a contradiction.
\end{proof}

The analogue of Proposition \ref{prop7} for $lr$-building blocks is obtained by duality, and the same goes for the compactness argument extending it from
gap pairs to gaps. We skip the details.
Sufficient conditions for the three-point
inequality will be discussed in Section \ref{sect_sufficient}.

\section{Convergence}
\label{sect_convergence_rl}
We start with a global convergence result.

\begin{theorem}
\label{main_convergence}
{\bf (Global convergence)}.
Let $a_k,b_k$ be a Bregman alternating sequence with gap $(K^*,r^*)$. 
Suppose the $rl$-angle condition and $rl$-three-point inequality are satisfied at every $(b^*,a^*)\in K^*$.
Then the sequence $b_k\in \cev{P}_B\circ \vec{P}_A(b_{k-1})$ converges.
\end{theorem}

\begin{proof}
1) Since $D(b_{k-1},a_{k-1}) \leq D(b_{k-1},a_{k}) \leq D(b_k,a_k)$, we have $D(b_k,a_k) \to \frac{1}{2}r^{*2}$ and also
$D(b_{k-1},a_k) \to \frac{1}{2}r^{*2}$ by monotone convergence.  Due to boundedness of the $a_k,b_k$ the set
$K^*$ of accumulation point pairs $(b^*,a^*)$ of the alternating sequence satisfying $b^*\sim a^*$  is compact, so
by the usual compactness argument
we find a neighborhood $W$ of $K^*$ on which the angle condition holds with the same $\gamma, \eta$ and de-singularizing function 
$\phi$, respectively, the same angle shrinking function $\sigma(s)=\phi'(s)^{-2}$ or  $\sigma(s)=\phi'(s)^{-2}s^{-1}$, simultaneously  for all $lr$-building blocks in $W$. Shrinking  $W$ further,  
and using Proposition \ref{prop7}, we may in addition arrange that  the three point inequality is satisfied with $\ell$ for all $lr$-building blocks
in $W$.
Since there are only finitely many iterates outside $W$, we may without loss assume
\begin{equation}
\label{above}
\phi'(\cev{D}_B(a_k)-\textstyle\frac{1}{2}r^{*2})^2 \cev{D}_B(a_k)  (1-\cos \alpha_k) \geq \gamma
\end{equation}
for all $k$, 
where $\alpha_k = \angle(b_{k-1}-a_k,b_k-a_k)$.
With the same argument    we can also assume that
the three-point inequality holds for  all $k$ with the same constant $\ell\in (0,1)$, hence we also have the four point inequality
\begin{equation}
\label{four_point_gap}
D(b_{k},a_{k}) \geq D(b_{k+1},a_{k+1}) + \ell D(b_k,b_{k+1}).
\end{equation}

2)
Due to our standing assumption $A,B \subset G$, we may further
assume that (\ref{bounds}) holds on $W$ with suitable constants $m,M$. 
Now
\begin{align}
\begin{split}
m^{-2} D(&b_{k-1},b_k)  \geq  \|b_{k-1}-b_k\|^2 \\
&=  \|b_{k-1}-a_k\|^2 + \|a_k-b_k\|^2 - 2\|b_{k-1}-a_k\|\|a_k-b_k\| \cos \alpha_k \\
&=  \left( \|b_{k-1}-a_k\|- \|a_k-b_k\|  \right)^2 + 2\|b_{k-1}-a_k\| \|a_k-b_k\| (1-\cos\alpha_k) \\
&\geq 2\|b_{k-1}-a_k\| \|a_k-b_k\| (1-\cos\alpha_k) \\
&\geq 2 M^{-1} D(b_{k-1},a_k)^{1/2} \|a_k-b_k\| (1-\cos \alpha_k)  \\
&\geq 2\ M^{-1} D(b_k,a_k)^{1/2} \|a_k-b_k\| (1-\cos\alpha_k) \\
&\geq 2M^{-2} D(b_k,a_k) (1-\cos\alpha_k)  \\
&\geq 2M^{-2} \gamma \phi'(\cev{D}_B(a_k) - \textstyle \frac{1}{2}r^{*2})^{-2}.
\end{split}
\end{align}
Here lines 1,5,7 use (\ref{bounds}),  
lines 2-4 concern the cosine theorem,
line 6 uses $D(b_k,a_k) \leq D(b_{k-1},a_k)$,  and the last line
uses the angle condition (\ref{above}).
We re-write this as
\begin{equation}
\label{need_den}
D(b_{k-1},b_k)^{1/2} \geq \frac{m\sqrt{2\gamma}}{M}  \phi'(\cev{D}_B(a_k) - \textstyle \frac{1}{2}r^{*2})^{-1}.
\end{equation}

3) Concavity of the de-singularizing function $\phi$ now implies
\begin{align*}
\phi(\cev{D}_B(a_k)-\textstyle\frac{1}{2}r^{*2})  - \phi(\cev{D}_B(a_{k+1})-\frac{1}{2}r^{*2})& \geq \phi'(\cev{D}_B(a_k) - \textstyle \frac{1}{2}r^{*2})
\left[\cev{D}_B(a_k)-\frac{1}{2}r^{*2} - (\cev{D}_B(a_{k+1})-\frac{1}{2}r^{*2})\right] \\
&=  \phi'(\cev{D}_B(a_k) - \textstyle \frac{1}{2}r^{*2}) \left[D(b_k,a_k) - D(b_{k+1},a_{k+1}) \right]\\
&\geq  \phi'(\cev{D}_B(a_k) - \textstyle \frac{1}{2}r^{*2})  \ell D(b_k,b_{k+1}) \\
&\geq m\sqrt{2\gamma} M^{-1} \ell \frac{D(b_{k},b_{k+1})}{D(b_{k-1},b_k)^{1/2}}.
\end{align*} 
Here the third line uses the four point inequality (\ref{four_point_gap}), while the last line uses (\ref{need_den}). 
Setting $C = M/m\sqrt{2\gamma}\ell$, we have
$$
C \left[ \phi(\cev{D}_B(a_k)-\textstyle\frac{1}{2}r^{*2}) - \phi(\cev{D}_B(a_{k+1})-\frac{1}{2}r^{*2}) \right] D(b_{k-1},b_k)^{1/2} \geq D(b_k,b_{k+1}).
$$
Since $a^2\leq bc$ implies $a\leq \frac{1}{2}b + \frac{1}{2}c$ for positive $a,b,c$, we deduce
\begin{equation}
\label{step}
D(b_{k},b_{k+1})^{1/2} \leq \frac{1}{2} D(b_{k-1},b_k)^{1/2} + \frac{C}{2}  \left[ \phi(\cev{D}_B(a_k)-\textstyle\frac{1}{2}r^{*2}) - \phi(\cev{D}_B(a_{k+1})-\frac{1}{2}r^{*2}) \right].
\end{equation}
Summing this from $k=1$ to $n$ gives
$$
\sum_{k=1}^n D(b_k,b_{k+1})^{1/2} \leq \frac{1}{2} \sum_{k=1}^n D(b_{k-1},b_k)^{1/2} + \frac{C}{2} \left[ \phi(\cev{D}_B(a_1)-\textstyle\frac{1}{2}r^{*2}) - \phi(\cev{D}_B(a_{n+1})-\frac{1}{2}r^{*2}) \right].
$$
Re-arranging, and multiplying by $2$, we obtain
\begin{align*}
\sum_{k=1}^n D(b_k,b_{k+1})^{1/2} & \leq D(b_0,b_1)^{1/2}
+ C \left[ \phi(\cev{D}_B(a_1)-\textstyle\frac{1}{2}r^{*2}) - \phi(\cev{D}_B(a_{n+1})-\frac{1}{2}r^{*2}) \right] - D(b_{n},b_{n+1})^{1/2}\\
&\leq D(b_0,b_1)^{1/2} + C \phi(\cev{D}_B(a_1)-\textstyle\frac{1}{2}r^{*2}) .
\end{align*}
Using (\ref{bounds}), this implies
$$
\sum_{k=1}^n \|b_{k}-b_{k+1}\| \leq m^{-1}M \|b_0-b_1\| + m^{-1}  C \phi(\cev{D}_B(a_1)-\textstyle\frac{1}{2}r^{*2}) .
$$
This proves convergence of the series $\sum_{k=1}^\infty \|b_k-b_{k+1}\|$, hence the sequence $b_k$ is Cauchy and converges to some $b^*\in B$.
\end{proof}

\begin{remark}
For $r^* > 0$ the proof does not assure convergence of the $a_k$, but
all accumulation points 
$a^*$ of the $a_k$ satisfy $D(b^*,a^*) = \frac{1}{2}r^{*2}$ and $b^* \sim a^*$, i.e., lie on the boundary of $\vec{\mathcal B}(b^*,r^*)$, and the gap has the form
$K^* = \{b^*\} \times A^*$. 
In the feasible case $r^*=0$, convergence of the $a_k$ is guaranteed, but here we have an even stronger result:
\end{remark}

\begin{theorem}
\label{thm_attract}
{\bf (Convergence by attraction)}.
Let $x^*\in A \cap B$, and suppose rl-angle condition and rl-three point estimate are satisfied at $x^*$. Then there exists a neighborhood $V$
of $x^*$ such that every Bregman alternating sequence which enters $V$ converges to some point in the intersection.
\end{theorem}

\begin{proof}
1) 
By Proposition \ref{prop_gap_angle} the angle condition (\ref{angle}) holds on a neighborhood $U$ of  $x^*\in G$ with the 
shrinking function $\sigma(s) = 1/s\phi'(s)^2$ and a constant $\gamma$. In addition, 
$U$ may be chosen such that every building block
$b \stackrel{r}{\longrightarrow} a^+ \stackrel{l}{\longrightarrow} b^+$ in $U$ satisfies the three point estimate
(\ref{rlr_3_pt}) with $0 < \ell < 1$. From the three-point inequality we immediately obtain the  following four-point-inequality
\begin{equation}
\label{four_point}
D(b,a) \geq D(b^+,a^+) + \ell D(b,b^+)
\end{equation}
for building blocks $a \stackrel{l}{\longrightarrow} b \stackrel{r}{\longrightarrow} a^+ \stackrel{l}{\longrightarrow} b^+$
with $b,a^+,b^+\in U$.

Since the alternating sequence including accumulation points is contained in the interior of dom$f$,  we may assume
that (\ref{bounds}) with constants $m,M$  is satisfied in a neighborhood of the set of iterates. In particular, we assume
that it is satisfied on the neighborhood $U$
of $x^*$.  
Let $U = B(x^*,\epsilon)$ without loss.  
Now define 
$$
C =   M/m \sqrt{2\gamma} \ell.    
$$

2)
Choose $\delta > 0$ such that the following are satisfied:
\begin{align}
\label{38}
\begin{split}
&(m^{-5}M^5+m^{-4}M^4  +  m^{-3}M^3  + m^{-2}M^2 +m^{-1} M +1) \delta < \epsilon/2\\
& (m^{-1} + m^{-2}M + m^{-3}M^2+m^{-4}M^3)C \phi(\xi) < \epsilon/2
\end{split}
\end{align}
for all $|\xi| < M^2\delta^2$.
The latter is possible due to $\phi(0)=0$ and continuity of $\phi$ at $0$.
Now let $V = B(x^*,\delta)$. We claim that if the alternating sequence 
$a_{k} \stackrel{l}{\longrightarrow} b_{k} \stackrel{r}{\longrightarrow} a_{k+1} \stackrel{l}{\longrightarrow} b_{k+1}$
enters $V$, then
it converges to a point $b^\sharp\in A \cap B$. Relabeling the sequence,
we may  assume that $b_0 \in V$. The case where the $a_k$ reach $V$ first is treated analogously.

3)
We shall prove by induction that for every $k\geq 1$,
\begin{equation}
\label{induction1}
b_0,a_1,b_1,\dots, a_k,b_{k},a_{k+1},b_{k+1} \in U
\end{equation}
and
\begin{equation}
\label{induction2}
\sum_{j=1}^{k} D(b_{j},b_{j+1})^{1/2} \leq \frac{1}{2} \sum_{j=1}^{k} D(b_{j-1},b_j)^{1/2} + \frac{C}{2}  \left[ \phi(\cev{D}_B(a_1)) - \phi(\cev{D}_B(a_{k+1})) \right].
\end{equation}

Let us first prove (\ref{induction1}) for $k=1$. This means we have to show $a_1,b_1,a_2,b_2\in U$.
We have $D(b_0,a_1) \leq D(b_0,x^*)$ due to $x^*\in A$, hence
$m\|b_0-a_1\| \leq D(b_0,a_1)^{1/2} \leq D(b_0,x^*)^{1/2} \leq M \|b_0-x^*\|$, giving $\|b_0-a_1\| \leq m^{-1}M\delta$. Then
$\|a_1-x^*\| \leq \|a_1 - b_0\| + \|b_0-x^*\| \leq (m^{-1} M +1) \delta < \epsilon$ using (\ref{38}), which is the first claim.

Now $D(b_1,a_1) \leq D(b_0,a_1)$, hence
$m\|b_1-a_1\| \leq M \|b_0-a_1\| \leq m^{-1} M^2\delta$. Then $\|b_1-a_1\| \leq M^2m^{-2} \delta$, giving
$\|b_1-x^*\|\leq \|b_1-a_1\| + \|a_1-x^*\| \leq  (m^{-2}M^2 +m^{-1} M+1) \delta < \epsilon$, again using (\ref{38}). This is the second statement in $(\ref{induction1})_1$.

Next $D(b_1,a_2) \leq D(b_1,a_1)$, which gives
$m\|b_1-a_2\|\leq M\|b_1-a_1\| \leq m^{-2}M^3 \delta$, hence $\|b_1-a_2\|\leq m^{-3}M^3 \delta$.
Then $\|a_2-x^*\| \leq \|a_2-b_1\| + \|b_1-x^*\| \leq (m^{-3}M^3  + m^{-2}M^2 +m^{-1} M + 1)\delta   < \epsilon$ from (\ref{38}), which is the third
statement in $(\ref{induction1})_1$.

Finally, from $D(b_2,a_2) \leq D(b_2,a_1)$
we get $m\|b_2-a_2\| \leq M \|b_2-a_1\|$, hence $\|b_2-a_2\|\leq m^{-4}M^4 \delta$, so that
$\|b_2-x^*\| \leq \|b_2-a_2\| + \|a_2-x^*\| < (m^{-4}M^4  +  m^{-3}M^3  + m^{-2}M^2 +m^{-1} M +1) \delta < \epsilon$, once again via (\ref{38}).
That proves $a_1,b_1,a_2,b_2\in U$.

4)
Before proving $(\ref{induction2})_1$, let us first
do the induction step.
Suppose $(\ref{induction1})_{k-1}$ and $(\ref{induction2})_{k-1}$ are satisfied.
We have to prove $(\ref{induction1})_k$ and $(\ref{induction2})_k$.
We first check (\ref{induction1}) at $k$. By $(\ref{induction1})_{k-1}$
we know that $b_0,a_1,b_1,\dots,a_k,b_k\in U$, so it remains to
prove
$a_{k+1},b_{k+1}\in U$. Now observe that $(\ref{induction2})_{k-1}$ implies
\begin{align*}
\sum_{j=1}^{k-1} D(b_j,b_{j+1})^{1/2} &\leq D(b_0,b_1)^{1/2} +  {C}  \left[ \phi(\cev{D}_B(a_1)) - \phi(\cev{D}_B(a_{k})) \right] - D(b_{k-1},b_k)^{1/2}\\
&\leq D(b_0,b_1)^{1/2} +  {C}  \phi(\cev{D}_B(a_1)).
\end{align*}
Using (\ref{bounds})  
this implies
$$
\sum_{j=1}^{k-1} \|b_{j-1}-b_j\|\leq m^{-1} M \|b_0-b_1\| +m^{-1} C  \phi(\cev{D}_B(a_1)).
$$
Using this, we have
\begin{align*}
\|b_k-x^*\| &\leq \|b_k-b_1\| + \|b_1 - x^*\| \leq \sum_{j=1}^{k-1} \|b_j-b_{j+1}\| + \|b_1-x^*\| \\
&\leq m^{-1}M \|b_0-b_1\| + m^{-1}C  \phi(\cev{D}_B(a_1)) + \|b_1-x^*\| \\
&\leq m^{-1}M \|b_0-a_1\| + m^{-1}M \|a_1-b_1\| +  m^{-1}C  \phi(\cev{D}_B(a_1)) + (m^{-2}M^2 + m^{-1}M + 1)\delta\\
& < m^{-2}M^2 \delta + m^{-3}M^3 \delta + m^{-1}C  \phi(\cev{D}_B(a_1)) + (m^{-2}M^2 + m^{-1}M + 1)\delta\\
&= (m^{-3} M^3 + 2m^{-2}M^2 + m^{-1} M + 1)\delta + m^{-1}C  \phi(\cev{D}_B(a_1)).
\end{align*}
Now $D(b_k,a_{k+1}) \leq D(b_k,x^*)$, hence  
$\|b_k-a_{k+1}\| \leq m^{-1}M \|b_k-x^*\| 
\leq (m^{-4}M^4 +\dots +m^{-1}M)\delta + m^{-2}MC \phi(\cev{D}_B(a_1))$.
Then $\|a_{k+1}-x^*\| \leq \|a_{k+1}-b_k\| + \|b_k-x^*\| \leq 
(m^{-4}M^4+\dots+1) \delta + (m^{-1}+m^{-2}M)C\phi(\cev{D}_B(a_1))
< \epsilon/2 + \epsilon/2 = \epsilon$, where we use (\ref{38}), being allowed to do so due to $\cev{D}_B(a_1) < M^2\delta^2$.
Namely,
$\cev{D}_B(a_1) = D(b_1,a_1) \leq D(b_0,a_1) \leq D(b_0,x^*) \leq M^2\|b_0-x^*\|^2 \leq M^2\delta^2$,  which bounds the term
$(m^{-1}+m^{-2}M)C  \phi(\cev{D}_B(a_1))$ according to the second condition in (\ref{38}).

Finally, $D(b_{k+1},a_{k+1}) \leq D(b_k,a_{k+1})$, so
$\|b_{k+1}-a_{k+1}\|\leq m^{-1}M \|b_k-a_{k+1}\| \leq (m^{-5}M^5+\dots+m^{-1}M)\delta + m^{-3}M^2C\phi(\cev{D}_B(a_1))$, which gives
$\|b_{k+1}-x^*\| \leq \|b_{k+1}-a_{k+1} \| + \|a_{k+1}-x^*\| \leq
(m^{-5}M^5+\dots+1)\delta + (m^{-1} + m^{-2}M + m^{-3}M^2)C \phi(\cev{D}_B(a_1)) 
< \epsilon/2+\epsilon/2 = \epsilon$, using both estimates in (\ref{38}).
That proves $(\ref{induction1})_k$.
\\

6)
Now let us prove $(\ref{induction2})_k$.
From the angle condition (\ref{angle}) with $\sigma(s)=\phi'(s)^{-2}s^{-1}$ we have
\begin{equation}
\label{angle_k}
 \phi'(\cev{D}_B(a_k))^2 \cev{D}_B(a_k) (1-\cos\alpha_k)\geq \gamma,
\end{equation}
where $\alpha_k = \angle(b_{k-1}-a_k,b_k-a_k)$. 
We also have the four-point estimate
\begin{equation}
\label{four_k}
D(b_{k+1},a_{k+1}) + \ell D(b_k,b_{k+1}) \leq D(b_k,a_k)
\end{equation}
because $D(b_k,a_{k+1}) \leq D(b_k,a_k)$.  Now we 
invoke precisely the same estimation as used in the proof of Theorem \ref{main_convergence},
which via (\ref{above}) and (\ref{four_point_gap}) led to (\ref{need_den}). Here we use instead (\ref{angle_k})
and (\ref{four_k}) to derive (using $r^*=0$):
\begin{equation}
\label{21}
D(b_{k-1},b_k)^{1/2} \geq \frac{m \sqrt{2\gamma}}{M} \phi'(\cev{D}_B(a_k))^{-1}.
\end{equation}

Now by concavity of $\phi$, the four point estimate (\ref{four_k}), and (\ref{21}), we get
\begin{align*}
\phi(\cev{D}_B(a_k)) - \phi(\cev{D}_B(a_{k+1})) &\geq \phi'(\cev{D}_B(a_k)) \left( D(b_k,a_k) - D(b_{k+1},a_{k+1} \right)\\
&\geq  \phi'(\cev{D}_B(a_k))  \ell D(b_k,b_{k+1})\\
&\geq m \sqrt{2\gamma} M^{-1} \ell \frac{D(b_k,b_{k+1})}{D(b_{k-1},b_k)^{1/2}}.
\end{align*}
With  $C = M/m \sqrt{2\gamma} \ell$  this becomes
$$
C \left[ \phi(\cev{D}_B(a_k)) - \phi(\cev{D}_B(a_{k+1})) \right] D(b_{k-1},b_k)^{1/2} \geq D(b_k,b_{k+1}).
$$
Since $a^2\leq bc$ implies $a\leq \frac{1}{2}b + \frac{1}{2}c$ for positive $a,b,c$, we deduce
\begin{equation}
\label{before}
D(b_k,b_{k+1})^{1/2} \leq  \frac{1}{2} D(b_{k-1},b_k)^{1/2} + \frac{C}{2}  \left[ \phi(\cev{D}_B(a_k)) - \phi(\cev{D}_B(a_{k+1})) \right].
\end{equation}
By the induction hypothesis $(\ref{induction2})_{k-1}$ we have
$$
\sum_{j=1}^{k-1} D(b_{j},b_{j+1})^{1/2} \leq \frac{1}{2} \sum_{j=1}^{k-1} D(b_{j-1},b_j) ^{1/2}+ \frac{C}{2}  \left[ \phi(\cev{D}_B(a_1)) - \phi(\cev{D}_B(a_{k})) \right].
$$
Adding this and (\ref{before}) gives  $(\ref{induction2})_k$ at stage $k$.

7) It remains to prove (\ref{induction2}) at $k=1$. This can be done by following the same steps as in 6) with $k=1$.

8)
Having proved $(\ref{induction2})_k$ for all $k$, we see that the series $\sum_{k=1}^\infty D(b_{k-1},b_k)^{1/2}$ converges, and that all iterates
stay in $U$. Hence via (\ref{bounds}) the series $\sum_{k=1}^\infty \|b_{k-1}-b_k\|$ converges as well, hence the $b_k$ form a Cauchy sequence, which converges
to some $b^\sharp\in B \cap U$.

9)
Convergence of the $a_k$ is obtained as follows. We have $D(b_k,a_{k+1}) \leq D(b_k,b^\sharp)$ due to $b^\sharp\in A \cap B$ and $a_{k+1}\in \vec{P}_A(b_k)$, hence
$\|a_{k+1}-b_k\|\leq m^{-1}M \|b_k-b^\sharp\|$, and then $\|a_{k+1}-b^\sharp\|\leq \|a_{k+1}-b_k\| + \|b_k-b^\sharp\| \leq (1+m^{-1}M)\|b_k-b^\sharp\|$.
\end{proof}

\begin{corollary}
\label{cor2}
Suppose $A,B,f$ are definable and $B$ is prox-regular at $x^*\in A \cap B$. Then
there exists a neighborhood $V$ of $x^*$ such that
every Bregman alternating sequence $a_{k} \stackrel{l}{\longrightarrow} b_{k} \stackrel{r}{\longrightarrow} a_{k+1}$
which reaches $V$ converges to a point $b^\sharp\in A \cap B$.
\end{corollary}

We pick up the thread of local projections from Section \ref{sect_local}.

\begin{corollary}
\label{cor_local}
Consider a local Bregman alternating sequence 
$b_{k-1} \stackrel{\,r\,}{\dashrightarrow} a_k \stackrel{l}{\longrightarrow} b_k$
with gap $K^*$. Suppose $A,B,f$ are definable.  Let the $rl$-three point inequality be satisfied
on $K^*$. Then the sequence $b_k$ converges.
\end{corollary}

\begin{proof}
Since $a_k,b_k$ is a standard Bregman alternating sequence between the sets $A^s,B^s$,
we know from Section \ref{sect_local} that the
angle condition, which holds throughout $A,B$, remains in place between $A^s$ and $B^s$. Let $K^*$ be the gap generated by the $a_k,b_k$, then 
$K^*$ is also the gap of $a_k,b_k$ when the latter is considered alternating between  $A^s$ and $B^s$.
By hypothesis the $rl$-three point inequality holds on $K^*$, and since it is also expressed in terms of building blocks, it
remains true for $a_k,b_k$ alternating between $A^s,B^s$. Altogether, by Theorem \ref{main_convergence}, the sequence
$b_k$ converges.
\end{proof}

This is important from a practical point of view when $A$ is not convex, as we then may want to solve (\ref{right}) with a local NLP-solver, starting
with the last $a$ as initial guess, thereby assuring descent (\ref{decrease}).

\section{dual convergence}
\label{sect_dual_convergence}
We continue to consider the alternating sequence under the form
$$
a \stackrel{l}{\longrightarrow} b \stackrel{r}{\longrightarrow} a^+ \stackrel{l}{\longrightarrow} b^+
$$
where 
$D(b^+,a^+)\leq D(b,a^+) \leq D(b,a)$. However, we now
reverse the roles of $A$ and $B$, i.e.,  we assume that the dual angle condition (\ref{angle_lr})
and the dual
three point estimate (\ref{lrl_3_pt})
are satisfied, now for $lr$-building blocks $a \stackrel{l}{\longrightarrow} b \stackrel{r}{\longrightarrow} a^+$.

\begin{theorem}
\label{theorem3}
Let $a_k,b_k$ be a Bregman alternating sequence with gap $(K^*,r^*)$. Suppose the $lr$-angle condition and $lr$-three point inequality are satisfied at every
pair $(b^*,a^*)\in K^*$. Then the sequence $a_k = \vec{P}_A\circ \cev{P}_B(a_{k-1})$ converges.
\end{theorem}

\begin{proof}
We use duality to obtain the mirror sequence
$a_k^*,b_k^*$ of the $a_k,b_k$. Then $(K^*,r^*)$ is mapped into the dual gap $(\nabla f(K^*),r^*)$,
and by amenability of the angle condition and the three-point inequality,
the dual sequence now satisfies the $rl$-angle condition and $rl$-three-point inequality at $\nabla f(K^*)$. Therefore
$b_k^*\in\mbox{${\cev{P}}^{{}^{{}^*}}_{B^*}$}\circ$$\mbox{${\vec{P}}^{{}^{{}^*}}_{A^*}(b^*_{k-1})$}$ converges by Theorem \ref{main_convergence}.
Mapping this back under $\nabla f^*$ yields convergence of $a_k = \vec{P}_A\circ \cev{P}_B(a_{k-1})$.
\end{proof}

It is again possible to obtain convergence by attraction using duality.

\begin{corollary}
Let $\bar{x}\in A \cap B$ and suppose $lr$-angle condition and $lr$-three-point inequality
are satisfied at $\bar{x}$. Then there exists a neighborhood $V$ of $\bar{x}$ such that every  Bregman alternating sequence which enters $V$ 
converges to some point in the intersection. 
\end{corollary}

Naturally, we can also address the case of local projections 
$a \stackrel{\,l\,}{\dashrightarrow} b \stackrel{r}{\longrightarrow} a^+$ via duality. We skip the details.

\section{Speed of convergence}
\label{sect_speed}
We consider the case of the \L ojasiewicz inequality, where the de-singularizing function is
$\phi(s) = s^{1-\theta}$ for some $\theta\in  [\frac{1}{2},1)$. In that case,
worst case convergence rates can be obtained. 

\begin{corollary}
\label{cor_speed}
Under the hypotheses of Theorem {\rm \ref{main_convergence}}, suppose the de-singularizing function is of the form
$\phi'(s) = s^{-\theta}$ for some $\theta\in (\frac{1}{2},1)$. Then the speed of convergence of the sequence
$b_k=\cev{P}_B\circ \vec{P}_A(b_{k-1})$ is $\|b_k-b^*\| = O(k^{-\rho})$ with
$\rho = \frac{1-\theta}{2\theta-1}\in (0,\infty)$. 
When $\theta=\frac{1}{2}$ the speed is R-linear. In the feasible case, the $a_k$ converge to $b^*\in A\cap B$ with the same speed.
\end{corollary}

\begin{proof}
Summing (\ref{step}) from $k=N$ to $k=K$ gives
\begin{align}
\begin{split}
-\frac{1}{2} D(b_{N-1},b_N)^{1/2} + &\frac{1}{2} \sum_{k=N}^{K-1} D(b_k,b_{k+1})^{1/2} + D(b_K,b_{K+1})^{1/2}
\\
&\leq \frac{C}{2} \left[  \phi(\cev{D}_B(a_N)-\textstyle\frac{1}{2}r^{*2}) - \phi(\cev{D}_B(a_{K+1})-\textstyle\frac{1}{2}r^{*2}) \right].
\end{split}
\end{align}
Passing to the limit $K\to \infty$ gives
$$
-\frac{1}{2} D(b_{N-1},b_N)^{1/2} + \frac{1}{2} \sum_{k=N}^{\infty} D(b_k,b_{k+1})^{1/2}
\leq \frac{C}{2} \phi(\cev{D}_B(a_N)-\textstyle\frac{1}{2}r^{*2}). 
$$
Introducing $S_N=\sum_{k=N}^\infty D(b_k,b_{k+1})^{1/2}$, this reads
$$
-\frac{1}{2}(S_{N-1}-S_N) + \frac{1}{2} S_N \leq \frac{C}{2}  \phi(\cev{D}_B(a_N)-\textstyle\frac{1}{2}r^{*2}).
$$
On the other hand, (\ref{need_den}) gives
$$
\phi'(\cev{D}_B(a_N)-\textstyle\frac{1}{2}r^{*2})^{-1} \leq \displaystyle\frac{M}{m\sqrt{2\gamma}} D(b_{N-1},b_N)^{1/2} = \frac{M}{m\sqrt{2\gamma}} (S_{N-1}-S_N).
$$
Now by hypothesis we have $\phi'(s) = s^{-\theta}$ for $\theta\in [\frac{1}{2},1)$, hence $\phi(s) = (1-\theta)^{-1} s^{1-\theta}$. Therefore
$\left[ \phi'(s)^{-1} \right]^{\frac{1-\theta}{\theta}}=s^{1-\theta}= (1-\theta) \phi(s)$.
Hence $\phi(\cev{D}_B(a_N)-\textstyle\frac{1}{2}r^{*2})
\leq(1-\theta)^{-1}  [\phi'(\cev{D}_B(a_N)-\textstyle\frac{1}{2}r^{*2})^{-1}]^{\frac{1-\theta}{\theta}}  \leq (1-\theta)^{-1} \left( \frac{M}{m\sqrt{2\gamma}} \right)^{\frac{1-\theta}{\theta}} (S_{N-1}-S_N)^{\frac{1-\theta}{\theta}}$.
Substituting this gives
\begin{equation}
\label{bifork}
\frac{1}{2}S_N \leq C'(S_{N-1}-S_N)^{\frac{1-\theta}{\theta}} + \frac{1}{2} (S_{N-1}-S_N)
\end{equation}
with $C' = \frac{C}{2} (1-\theta)^{-1} \left( \frac{M}{m\sqrt{2\gamma}} \right)^{\frac{1-\theta}{\theta}}$.
Now for $\theta > \frac{1}{2}$ we have $\frac{1-\theta}{\theta} < 1$, so that the first term on the right hand side
dominates the second term. Therefore there exists another constant $C''$ such that
$$
S_N^{\frac{\theta}{1-\theta}} \leq C''(S_{N-1}-S_N).
$$
From here we follow precisely the argument in \cite[Cor. 4 (24) ff]{noll}, where this leads to an estimate of the form
$$
S_N \leq C''' N^{-\frac{1-\theta}{2\theta-1}}
$$
for another constant $C'''$. Using (\ref{bounds}), this shows $\widetilde{S}_N :=\sum_{k=N}^\infty \|b_k-b_{k+1}\| \leq m^{-1} S_N \leq m^{-1}C''' N^{-\frac{1-\theta}{2\theta-1}}$,
and since $\|b_N-b^*\| \leq \widetilde{S}_N$ by the triangle inequality,
we get the desired estimate $\|b_k-b^*\| = O(k^{-\rho})$ with $\rho = \frac{1-\theta}{2\theta-1}$. 

Now consider the case $\theta=\frac{1}{2}$, then (\ref{bifork}) turns into
$$
\frac{1}{2}
S_N \leq C' (S_{N-1}-S_N) + \frac{1}{2}(S_{N-1}-S_N),
$$
hence
$$
S_N \leq \frac{1+2C'}{2+2C'} S_{N-1},
$$
which gives Q-linear speed $S_N \to 0$, hence R-linear speed $\widetilde{S}_N \to 0$, and then also R-linear speed $\|b_N-b^*\|\to 0$.

Finally, in the feasible case we get the same speed of convergence for the $a_k$ from part 9) of the proof of Theorem \ref{thm_attract}.
\end{proof}

\begin{remark}
This may be compared to  \cite{kunstner}, where for a non-convex version of the EM algorithm for exponential families
a global estimate  $D(b_k,b_{k+1})^{1/2}  = O(k^{-1/2})$ is obtained without use of the K\L-inequality, with $D$ the Kullback-Leibler divergence.  
See also Example \ref{ex_2} in Section \ref{sect_examples}, and Section \ref{sect_em}.
\end{remark}

\begin{remark}
Naturally, via duality, the same speed of convergence is obtained for Theorem \ref{theorem3} if $\phi(s)=s^{1-\theta}$.
\end{remark}

\begin{corollary}
\label{cor_transversal}
{\bf (Linear convergence)}.
Suppose $B$ intersects $A$ $rl$-transversally at $\bar{x}\in A \cap B$, and the $rl$-three-point inequality is satisfied
at $\bar{x}$. Then there exists a neighborhood $V$ of $\bar{x}$ such that every Bregman alternating sequence which enters $V$
converges to some point in the intersection with R-linear speed. 
\end{corollary}

\begin{proof}
The neighborhood $V$ of $\bar{x}$ may be chosen such that the $rl$-three-point inequality holds on $V$.
By hypothesis we may also assure that the numerator $1-\cos\alpha$ in (\ref{angle}) stays bounded away from $0$ in $V$. Therefore we can allow a constant
as angle shrinking function $\sigma$. Now  $\sigma(s) = \phi'(s)^{-2} s^{-1}$ gives constant $\sigma$ as soon as $\phi'(s) \propto s^{-1/2}$, so the de-singularizing function
is $\phi(s) = s^{1/2}$. Then by Corollary \ref{cor_speed} convergence is R-linear near $\bar{x}$.
\end{proof}

The dual version of this result is also true.

\section{Sufficient conditions for the three-point inequality}
\label{sect_sufficient}
In this chapter we derive the three point inequality from conditions on the reach of the sets $A,B$.  Bregman reach had been introduced in Section \ref{sect_reach},
and extends the classical notion of reach \cite{federer}.  However, there are  other ways to extend the classical notion of reach to the Bregman setting, each
with advantages and inconveniences.

\subsection{Bregman reach larger than gap}
\label{sect_larger_1}
We start discussing left Bregman reach $\cev{R}(b^+,d)$, which we match with the distance
between $A$ and $B$ measured by $D(b^+,a^+)^{1/2}$. 

\begin{proposition}
\label{prop9}
Let $b^* \sim a^*$ with gap value $r^*\geq  0$ and suppose the left Bregman reach at $b^*\in B$ is at least $r > r^*$. 
Suppose $f$ is $1$-coercive.
Then there exist
$\delta > 0$ and  $0 < \ell < 1$ such that the three point inequality holds with $\ell$ for every building
block $b \stackrel{r}{\longrightarrow} a^+ \stackrel{l}{\longrightarrow}b^+$ with $(b^+,a^+)\in B(b^*,\delta) \times B(a^*,\delta)$.
\end{proposition}

\begin{proof}
1) Recall that the three point inequality holds trivially for every $\ell \in (0,1]$ if a building block satisfies
$\langle \nabla f(a^+)-\nabla f(b^+),b-b^+\rangle \leq 0$. We therefore assume $\langle \nabla f(a^+)-\nabla f(b^+),b-b^+\rangle > 0$ throughout.
Geometrically, this means that $b,a^+$ are strictly on the same side of the tangent hyperplane $H$ at $b^+$.

2)
We have $b^+\in \cev{P}_B(a^+)$ with $a^+\in A$, 
and the Bregman perpendicular $a_\lambda$ to $B$ at $b^+$ in direction $d=\nabla^2f(b^+)^{-1}(\nabla f(a^+)-\nabla f(b^+))$ satisfies 
\begin{equation}
\label{cos2}
\nabla f(a_\lambda) - \nabla f(b^+) = \lambda(\nabla f(a^+)-\nabla f(b^+)).
\end{equation}
Now consider
$b$ strictly on the same side of the tangent hyperplane as $a^+$. We claim that
there exists a parameter $\lambda(b) > 1$ for which the point $a_{\lambda(b)}$ on the geodesic
gives equality $D(b,a_{\lambda(b)}) = D(b^+,a_{\lambda(b)})$.
Indeed, the cosine theorem for Bregman distances in tandem with (\ref{cos2}) gives
\begin{equation}
\label{40}
D(b,a_\lambda) = D(b,b^+) + D(b^+,a_\lambda) - \lambda \langle \nabla f(a^+)-\nabla f(b^+),b-b^+\rangle,
\end{equation}
hence
\begin{align*}
\frac{D(b,a_\lambda) - D(b^+,a_\lambda)}{\lambda} &= \frac{D(b,b^+)}{\lambda} -  \langle \nabla f(a^+)-\nabla f(b^+),b-b^+\rangle\\
&\to - \langle \nabla f(a^+)-\nabla f(b^+),b-b^+\rangle<0
\end{align*}
as $\lambda \to \infty$, so that eventually $D(b^+,a_\lambda) > D(b,a_\lambda)$. Here we use the fact that due to $1$-coercivity of $f$ the  $a_\lambda$ are defined for all $\lambda \geq 0$, 
so that we may pass to the limit, and we use the fact that the limit term is negative.
Since at $a_\lambda|_{\lambda = 1}=a^+$ we have $D(b,a_1) > D(b^+,a_1)$,  the intermediate value theorem gives $\lambda=\lambda(b) \in (1,\infty)$ with equality.

3) 
Differentiating (\ref{cos2}) with respect to $\lambda$ gives
$\nabla^2f(a_\lambda) \frac{d}{d\lambda} a_\lambda = \nabla f(a^+)-\nabla f(b^+)$. Therefore
$\frac{d}{d\lambda} D(b^+,a_\lambda) = -\langle \frac{d}{d\lambda} a_\lambda, \nabla^2f(a_\lambda)(b^+-a_\lambda)\rangle
= -\langle \nabla^2f(a_\lambda)\frac{d}{d\lambda}a_\lambda,b^+-a_\lambda\rangle=\langle \nabla f(a^+)-\nabla f(b^+),a_\lambda - b^+\rangle>0$,
as all $a_\lambda$ are on the same side of the tangent hyperplane as $a^+$. This means $\lambda \mapsto D(b^+,a_\lambda)$ is strictly increasing,
hence the intermediate
value $\lambda(b)$ found above is unique.

4) From $D(b,a_{\lambda(b)}) = D(b^+,a_{\lambda(b)})$  we get
$D(b,b^+) = \lambda(b) \langle \nabla f(a^+)-\nabla f(b^+),b-b^+\rangle$ using (\ref{40}), hence the three point inequality
holds with $\ell(b) = 1-\lambda(b)^{-1}$ for every
building block
$b \stackrel{r}{\longrightarrow} a^+ \stackrel{l}{\longrightarrow}b^+$ with $b$ strictly on the same side of the hyperplane $H$ as $a^+$.

What remains to be shown is that there is one global $\ell$ which works for all these building blocks, or put differently, 
that the $\ell(b)=1-\lambda(b)^{-1}$ stay bounded away from $0$ even when their building blocks approach the gap.

From the construction we have $b,b^+\in \partial \cev{\mathcal B}(a_{\lambda(b)},r_{\lambda(b)})$,
so by the definition of the left reach
$r_{\lambda(b)} \geq \cev{R}(b^+,d) \geq r > r^*$  
for all such $(b^+,a^+) \in W$. 
Now assume contrary to what is claimed that there exist building blocks
$b_{k-1} \stackrel{r}{\longrightarrow} a_k \stackrel{l}{\longrightarrow}b_k$
with $a_k\to a^*$, $b_k\to b^*$, such that $\lambda(b_{k-1}) \to 1$. 
Then
$\nabla f(a_{\lambda(b_{k-1})}) = \nabla f(b_k) + \lambda({b_{k-1}}) (\nabla f(a_k)-\nabla f(b_k)) \to \nabla f(a^*)$.
Since $\nabla f$ is an diffeomorphism from int(dom $f$) onto int(dom$f^*$), this implies $a_{\lambda(b_{k-1})} \to a^*$.
Then $D(b_k,a_{\lambda(b_{k-1})}) \to D(b^*,a^*)=\frac{1}{2}r^{*2}$, hence $r_{\lambda(b_{k-1})} \to r^*$, a contradiction
with the above.
\end{proof}

\begin{remark}
In view of Proposition \ref{prop2}
it is unlikely that the result still holds without $1$-coercivity of $f$. This is why we consider alternative ways to define reach via $\widetilde{R}$ in the following sections.
\end{remark}

\subsection{Mobile reach larger than gap}
\label{sect_larger_2}
We use
a different notion of reach based on what we call a mobile euclidean norm. This  requires measuring the gap between the sets differently.
 
With every building block $b\stackrel{r}{\longrightarrow} a^+ \stackrel{l}{\longrightarrow} b^+$ we associate the euclidean norm
$\|x\|_{b^+}^2 = \langle x,\nabla^2f(b^+)x\rangle = \langle x,x\rangle_{b^+}$. 
Due to $\nabla^2f(b^+) \succeq \epsilon >  0$  for $b^+$ in a compact subset $K$ of the interior of dom$f$, we have an estimate of the form
\begin{equation}
\label{star_bounds}
m\|x\|_{b^+} \leq \|x\| \leq M\|x\|_{b^+} ,
\end{equation}
with $m,M$ depending only on $K$, and for every $b^+\in K$.  
The rationale of this norm stems from the fact that second-order Taylor-Young expansion of $f$ at $b^+$ gives
$$
D(b,b^+) = \textstyle\frac{1}{2}\|b-b^+\|_{b^+}^2 + o(\|b-b^+\|^2), \quad
D(b^+,a^+) = \frac{1}{2} \|b^+-a^+\|_{b^+}^2 + o(\|b^+-a^+\|^2).
$$
Here  the little-o terms may be made uniform on any compact set of $b^+$, see Section \ref{sect_uTJ}.

Let us fix some more terminology.

\begin{definition}
{\rm 
The proximal normal cone to $B$  with regard to $\|\cdot\|_{b^+}$ is $N_B^{p,b^+}$. The orthogonal projector with regard to
$\|\cdot\|_{b^+}$ is  $P_B^{b^+}$, the angle is $\angle_{b^+}$, and $\|\cdot\|_{b^+}$-balls are $B_{b^+}(x,r)$.
The reach of $B$ at $b^+\in B$ in direction $d\in N_B^{p,b^+}(b^+)\setminus\{0\}$ with regard to  $\|\cdot\|_{b^+}$ is 
$\widetilde{R}(b^+,d)$.
}
\end{definition}

This allows now the following

\begin{definition}
{\bf (Mobile reach)}.
\label{star_reach}
{\rm 
We say that $B$ has mobile reach at least $\widetilde{R}>0$ at $b^*\in B$, noted $\widetilde{R}(b^*) \geq \widetilde{R}$, 
 if there exists a neighborhood $U$ of $b^*$
such that for every $b^+\in P_B^{b^+}(a^+) \cap U$ for some $a^+\not \in B$
we have $\widetilde{R}(b^+,d) \geq \widetilde{R}$, where $d = \nabla^2f(b^+)^{-1}(\nabla f(a^+)-\nabla f(b^+))$. 
}
\end{definition}

\begin{remark}
As seen in Sections \ref{sect_rolling} and \ref{sect_reach},  $B$ has positive reach at $b^*\in B$ iff it has positive left Bregman reach at $b^*$.
Using (\ref{star_bounds}), this is now equivalent to having positive mobile reach. 
\end{remark}

\begin{proposition}
\label{prop10}
Let  $b^*\sim a^*$ with $\|\nabla f(b^*)-\nabla f(a^*)\|_{b^*} = \rho_*$. Suppose $B$ has mobile reach at least $\widetilde{R} > \rho_*$ at $b^*$.
Then there exist
$\delta > 0$ and $0 < \ell < 1$
such that every building block $b \stackrel{r}{\longrightarrow} a^+ \stackrel{l}{\longrightarrow}b^+$
with $(b^+,a^+)\in B(b^*,\delta)\times B(a^*,\delta)$  and $\|b-b^+\| < \delta$  satisfies the $rl$-three-point inequality with $\ell$.
\end{proposition}

\begin{proof}
Fix $\ell\in (0,1)$ such that $(1-\ell)^{-1} \rho_* < \widetilde{R}$. Then find $\epsilon > 0$
such that $(1+\epsilon)(1-\ell)^{-1} \rho_* < \widetilde{R}$.  Writing $\widetilde{R} = (1+\epsilon)(1-\ell)^{-1} \widetilde{\rho}$ therefore means $\rho_* < \widetilde{\rho}$.
Let $U$ be a neighborhood of $b^*$ as in Definition \ref{star_reach}. 
Shrink $U$
further  
until $(1+\epsilon) D(b,b^+) \geq \frac{1}{2}\|b-b^+\|_{b^+}^2$ for all $b,b^+\in U$, using uniform second order Taylor-Young expansion at $b^+$ on $U$
in tandem with (\ref{bounds}). 
Ready to shrink $U$ even further, combine it with
a neighborhood $V$ of $a^*$ such that $(b^+,a^+) \in U \times V$ implies $\|\nabla f(b^+) - \nabla f(a^+)\|_{b^+} <  \widetilde{\rho}$.  This is possible, because $\rho_* < \widetilde{\rho}$, and since the norm $\|\cdot\|_{b^+}$ depends continuously on $b^+$.

Now let $b \stackrel{r}{\longrightarrow} a^+ \stackrel{l}{\longrightarrow}b^+$ be a building block with $b,b^+\in U$ and $a^+\in V$, and put
$\bar{a} := b^++\nabla^2f(b^+)^{-1}(\nabla f(a^+)-\nabla f(b^+))$. Then 
\begin{align*}
\langle \nabla f(a^+)-\nabla f(b^+),b-b^+\rangle &= \langle \nabla^2f(b^+)^{-1}( \nabla f(a^+)-\nabla f(b^+)),b-b^+\rangle_{b^+}\\&=\langle \bar{a}-b^+,b-b^+\rangle_{b^+}\\
&= \|\bar{a}-b^+\|_{b^+} \|b-b^+\|_{b^+} \cos \beta,
\end{align*}
where $\beta = \angle_{b^+}(\bar{a}-b^+,b-b^+)$ is the angle in the euclidean geometry of $\|\cdot\|_{b^+}$. Now if $\cos\beta \leq 0$, the three-point inequality
is trivially satisfied, so we may assume $\cos\beta > 0$. Then $\beta < 90^\circ$, hence there exists a point
$\hat{a}$ on the proximal normal to $B$ at $b^+$ in the $\|\cdot\|_{b^+}$-geometry such that the triangle
$b,\hat{a},b^+$ is equilateral with two angles $\beta$ at the corners $b,b^+$, and two edges of the same $\|\cdot\|_{b^+}$-length $R$ joining $b,b^+$ to $\hat{a}$.  Hence $\hat{a} = b^++R d$,
where $d=(\bar{a}-b^+)/\|\bar{a}-b^+\|_{b^+}$
and $R = \|b-\hat{a}\|_{b^+}$.  Then $\|b-b^+\|_{b^+} = 2R\cos\beta$. But the ball $B_{b^+}(\hat{a},R)$ contains $b,b^+$, hence
$R \geq \widetilde{R}(b^+,d) \geq \widetilde{R} = (1+\epsilon)(1-\ell)^{-1}\widetilde{\rho}$. We deduce
$$
\|b-b^+\|_{b^+} \geq 2 (1+\epsilon)(1-\ell)^{-1} \widetilde{\rho} \cos\beta \geq 2(1+\epsilon) (1-\ell)^{-1} \|\nabla f(b^+)-\nabla f(a^+)\|_{b^+} \cos\beta.
$$
Hence
\begin{align*}
(1-\ell) D(b,b^+) & \geq (1-\ell)(1+\epsilon)^{-1} \textstyle\frac{1}{2} \|b-b^+\|_{b^+}^2\\
&\geq (1-\ell)(1+\epsilon)^{-1} (1+\epsilon) (1-\ell)^{-1} \|\nabla f(b^+)-\nabla f(a^+)\|_{b^+} \|b-b^+\|_{b^+}\cos\beta\\
&= \langle \nabla f(a^+)-\nabla f(b^+),b-b^+\rangle.
\end{align*}
That proves the $rl$-three point inequality.
\end{proof}

\begin{remark}
While Proposition \ref{prop9}  needs $1$-coercivity of $f$, which is not required here, we now need a nearness condition on $b,b^+$
to bring in the uniform Taylor-Young estimate, using that $f$ is of class $C^2$. This seems acceptable in view of Lemma \ref{lem7}. We will see the consequences right below.
\end{remark}

\begin{theorem}
{\bf (Global convergence)}.
Let $a_k,b_k$ be a Bregman alternating sequence with gap $(K^*,r^*)$. Suppose the $rl$-angle condition is satisfied at every $(b^*,a^*)\in K^*$.
Then the sequence $b_k \in \cev{P}_B \circ \vec{P}_A(b_{k-1})$ converges under any of the following conditions:
\begin{itemize}
\item[(i)] $B$ has left Bregman reach $\cev{R}(b^*) > r^*$ for every $b^* \sim a^*$ in $K^*$, and $f$ is $1$-coercive.
\item[(ii)] $B$ has mobile reach $\widetilde{R}(b^*) > \rho_*= \|\nabla f(b^*)-\nabla f(a^*)\|_{b^*}$ for all $b^*\sim a^*$ in $K^*$,  and the sequence satisfies $b_{k-1}-b_k\to 0$.
\end{itemize}
\end{theorem}

\begin{proof}
All we need  is assure the $rl$-three point inequality for  building blocks 
$b \stackrel{r}{\longrightarrow} a^+ \stackrel{l}{\longrightarrow} b^+$, as the rest is like in Theorem \ref{main_convergence}. We use Proposition \ref{prop9} for the left Bregman reach case (i), and
Proposition \ref{prop10} in case (ii). For the latter, by compactness we have $\widetilde{R} > \max\{\|\nabla f(b^*)-\nabla f(a^*)\|_{b^*}: b^* \sim a^*\}$, so that the three
point inequality holds for all gap pairs with the same $\ell$.
\end{proof}

\subsection{Slowly vanishing reach}
We have seen that positive reach at some $b^*\in B$, or prox-regularity of $B$ at $b^*$, could be expressed in four equivalent fashions, using $R,\cev{R}$, $\vec{R}$, and mobile reach $\widetilde{R}$.
Unfortunately, exact quantitative
relations among those four 
can only be obtained in the rough proportional sense of Section \ref{sect_rolling}. Better quantitative results can be 
obtained for zero gaps.  For those we may even allow $B$ to have {\it slowly vanishing reach} at $b^*\in B$, 
a notion introduced in \cite{noll,gerchberg} in the euclidean setting.  As this increases the chances to establish the three-point inequality, 
we adapt this to the Bregman context.
For an explanation of  what is meant by vanishing reach see    
also Example \ref{slow_reach}.
 
The following preparatory result conveys the fact that asymptotically as $r\to 0$, left Bregman balls $\cev{\mathcal B}(a^+,r)$ with $b^+\in \partial \cev{\mathcal B}(a^+,r)$
more and more resemble balls $B_{b^+}(\bar{a},\|\bar{a}-b^+\|_{b^+})$ in the euclidean geometry
$\langle x,y\rangle_{b^+} = \langle \nabla^2f(b^+)x,y\rangle$, where $\bar{a}= b^++\nabla^2f(b^+)^{-1}(\nabla f(a^+)-\nabla f(b^+))$.

\begin{lemma}
\label{lemma_den}
Let $x^*\in A \cap B$. Then
\begin{equation*}\liminf_{b^+\in \cev{P}_B(a^+), A \ni a^+\to x^*} \widetilde{R}(b^+,d) > 0 \mbox{ iff }
\liminf_{b^+\in \cev{P}_B(a^+), A \ni a^+\to x^*} \cev{R}(b^+,d)>0.\end{equation*}
When both tend to zero, we have
\begin{equation}\label{equal}\lim_{b^+\in \cev{P}_B(a^+), A \ni a^+\to x^*} \frac{\cev{R}(b^+,d)}{\widetilde{R}(b^+,d)}=1.\end{equation}
\end{lemma}

\begin{proof}
Note that $\widetilde{R}(b^+,d)$ stays away from 0 iff $\cev{R}(b^+,d)$ stays away from 0 by the results of Section \ref{sect_rolling} in tandem with (\ref{star_bounds}). 
Now suppose
both reach terms shrink to 0. 
Let $\cev{R}(b^+,d)$ with $d = \nabla^2f(b^+)^{-1}(\nabla f(a^+)-\nabla f(b^+))$ be realized at $\lambda > 1$ and radius $r_\lambda>0$, 
with $\cev{\mathcal B}(a_\lambda,r_\lambda)$ the largest left Bregman ball having $b^+$ on its boundary and no point of $B$ in its interior.

Working in the $\|\cdot\|_{b^+}$-geometry, uniform second-order Taylor-Young expansion of $f$ reads
$$
f(x+h) = f(x) + \langle \nabla^2f(b^+)^{-1}\nabla f(x),h\rangle_{b^+} + \textstyle\frac{1}{2} \langle h,\nabla^2f(b^+)^{-1} \nabla^2 f(x)h\rangle_{b^+} + o(\|h\|^2),
$$
hence the normal curvature of the left Bregman ball $\cev{\mathcal B}(a_\lambda,r_\lambda)$ at $x \in  \partial\cev{\mathcal B}(a_\lambda,r_\lambda)$ in unit tangential direction $v$ in the $\|\cdot\|_{b^+}$-geometry is
$$
\kappa_{n,b^+}(x,v) = \frac{\langle v,\nabla^2f(b^+)^{-1}\nabla^2f(x)v\rangle_{b^+}}{\| \nabla^2f(b^+)^{-1}(\nabla f(a_\lambda)-\nabla f(x))\|_{b^+}}.
$$
Taylor-Young expansion gives  $\nabla^2 f(b^+)^{-1}(\nabla f(a_\lambda)-\nabla f(x)) = a_\lambda - x + o(\|a_\lambda-b^+\|)+o(\|x-b^+\|)$. Since we assume $r_\lambda \to 0$ as the building block
approaches $x^*$, we have $a_\lambda \to x^*$ and $x \to x^*$ for $x\in \partial\cev{\mathcal B}(a_\lambda,r_\lambda)$. 
Hence on a sufficiently small neighborhood $U$ of $x^*$,
$(1-\epsilon) \|a_\lambda-x\|_{b^+} \leq \| \nabla^2f(b^+)^{-1}(\nabla f(a_\lambda)-\nabla f(x))\|_{b^+}
\leq (1+\epsilon) \|a_\lambda-x\|_{b^+}$ for $x,a_\lambda, b^+\in U$, using again that the little-o terms in the Taylor-Young expansion may be made uniform on a compact set of $b^+$ (Section \ref{sect_uTJ}). 
Shrinking $U$ further, we may arrange $\nabla^2f(b^+)^{-1} \nabla^2 f(x) = I_d + E$  with $\|E\|_{b^+}\leq \epsilon$ for 
$x,b^+\in U$. Then
$$
\frac{1-\epsilon}{(1+\epsilon) \|a_\lambda-x\|_{b^+}} \leq \kappa_{n,b^+}(x,v) \leq \frac{1+\epsilon}{(1-\epsilon) \|a_\lambda-x\|_{b^+}}
$$
using $\|v\|_{b^+}=1$.
This means the constants of Proposition \ref{prop_rolling}, applied in the $\|\cdot\|_{b^+}$-geometry,  are $\underline{c}=\frac{1-\epsilon}{1+\epsilon}$ and $\overline{c} = \frac{1+\epsilon}{1-\epsilon}$. 
Since the $\|\cdot\|_{b^+}$-euclidean ball with radius $\overline{c}r_\lambda$  contains the Bregman ball and touches it from outside at $b^+$, this ball contains also $b$, 
hence its radius is larger than the $\|\cdot\|_{b^+}$-reach at $b^+$:
$\frac{1+\epsilon}{1-\epsilon} \cev{R}(b^+,d)=\frac{1+\epsilon}{1-\epsilon} r_\lambda \geq \widetilde{R}(b^+,d)$.  On the other hand,
the smaller $\|\cdot\|_{b^+}$-ball with radius $\underline{c}r_\lambda$ is contained in the Bregman ball and touches it at $b^+$ from inside, hence contains no points of $B$
in its interior, whence $\frac{1-\epsilon}{1+\epsilon} r_\lambda \leq \widetilde{R}(b^+,d)$. 
This shows that numerator and  denominator in (\ref{equal}) agree asymptotically, which proves equality $1$ in the limit.
\end{proof}

\begin{definition}
{\bf (Slowly vanishing reach)}.
{\rm 
Let $x^*\in  A \cap B$. We say that $B$ has slowly vanishing reach with rate $\tau$ at $x^*$ with regard to $A$ if
\begin{equation}
\label{slowly_rl}
\tau := \limsup_{b^+\in \cev{P}_B(a^+), A \ni a^+\to x^*} \frac{D(b^+,a^+)^{1/2}}{\widetilde{R}(b^+,d)} < \frac{1}{\sqrt{2}}.
\end{equation}
}
\end{definition}

\noindent
Replacing $\widetilde{R}(b^+,d)$ by $\cev{R}(b^+,d)$ gives the same value $\tau$ by Lemma \ref{lemma_den}, hence an equivalent definition of slowly vanishing reach. 
We now see that $\widetilde{R}$
has the advantage over $\cev{R}$ that it gives the same information at points $x^*\in A \cap B$, while based on a euclidean norm, the inconvenience being
that it is a mobile one.

\subsection{Three point inequality from vanishing reach}
\label{sect_3_pt_suff}

In this section we show that the three-point inequality holds in the vicinity of a zero gap
$r^*=0$ even when we allow the reach
to shrink to zero, provided it shrinks slightly slower than the distance between the sets.

\begin{proposition}
\label{prop11}
Let $x^*\in A \cap B$, and suppose $B$ has slowly vanishing reach at $x^*$ with rate $\tau<{1}/{\sqrt{2}}$
with regard to $A$.
Then there exist a neighborhood $U$ of $x^*$
such that every building block
$b \stackrel{r}{\longrightarrow} a^+ \stackrel{l}{\longrightarrow}b^+$ with $b,a^+,b^+\in U$
satisfies the rl-three-point inequality with $\ell$ as long as $\ell$ satisfies  $\ell < 1 - \sqrt{2} \tau$.
\end{proposition}

\begin{proof}
1) 
The cosine theorem for Bregman distances gives
$$
D(b,a^+) = D(b^+,a^+) + D(b,b^+) - \langle \nabla f(a^+) - \nabla f(b^+),b-b^+\rangle.
$$
Therefore the rl-three-point inequality holds with $0 < \ell < 1$ iff
$$
(1-\ell)D(b,b^+) \geq \langle \nabla f(a^+)-\nabla f(b^+),b-b^+\rangle.
$$
This holds regardless of $\ell$
when $\langle \nabla f(a^+) - \nabla f(b^+),b-b^+\rangle \leq 0$.
We therefore
concentrate on the case  $\langle \nabla f(a^+) - \nabla f(b^+),b-b^+\rangle > 0$.

2) Let $\ell$ be as in the statement, and choose $\tau'$ such that $\tau < \tau' < \frac{1-\ell}{\sqrt{2}}$. 
Then choose $\epsilon > 0$ such that $\tau'(1+\epsilon)^4 < \frac{1-\ell}{\sqrt{2}}$. By the definition of $\tau$
we can  find a neighborhood $W$ of  $x^*$
such that 
\begin{equation}
\label{limsup_2}
\frac{D(b^+,a^+)^{1/2}}{\widetilde{R}(b^+,d)} \leq  \tau' 
\end{equation}
for all $a^+ \stackrel{r}{\longrightarrow} b^+$ with  $b^+,a^+\in W$.

We put
$\bar{a} = b^+ + \nabla^2f(b^+)^{-1}(\nabla f(a^+)-\nabla f(b^+))$. Shrinking $W$ further such that
$b,b^+$ are forced sufficiently close, using
uniform Taylor-Young expansion $D(b,b^+) = \frac{1}{2}\|b-b^+\|_{b^+}^2 + o(\|b-b^+\|^2)$ in tandem with (\ref{bounds}), 
we may arrange the following:
\begin{align}
\label{26_1}
\begin{split}
{\rm (a)} \quad & \textstyle\frac{1}{2} \|b-b^+\|_{b^+}^2 \leq (1+\epsilon)^2 D(b,b^+),\\
{\rm (b)} \quad & \textstyle \frac{1}{2} \|b^+-a^+\|_{b^+}^2 \leq (1+\epsilon)^2 D(b^+,a^+)  \\
{\rm (c)} \quad & \|b^+-\bar{a}\|_{b^+} \leq (1+\epsilon) \|b^+-a^+\|_{b^+}   \\
\end{split}
\end{align}
for all building blocks with $a^+,b^+,\bar{a},b\in W$.
Let the angle $\beta = \angle_{b^+}(b-b^+,\bar{a}-b^+)$ be
taken with regard to the $\|\cdot\|_{b^+}$-geometry. 
Then,
\begin{align}
\label{27_1}
\begin{split}
\langle \nabla f(a^+)-\nabla f(b^+),b-b^+\rangle &= \langle \nabla^2f(b^+)^{-1}(\nabla f(a^+)-\nabla f(b^+)),b-b^+\rangle_{b^+} \\
&= \|\nabla^2f(b^+)^{-1}(\nabla f(a^+)-\nabla f(b^+))\|_{b^+} \|b-b^+\|_{b^+} \cos\beta \\
&= \|\bar{a}-b^+\|_{b^+}\|b-b^+\|_{b^+} \cos \beta\\
&\leq (1+\epsilon) \|b^+-a^+\|_{b^+} \|b-b^+\|_{b^+} \cos\beta \\
&\leq (1+\epsilon)^2 \|b^+-a^+\|_{b^+} \sqrt{2} D(b,b^+)^{1/2} \cos\beta\\
&\leq (1+\epsilon)^3 \sqrt{2} D(b^+,a^+)^{1/2} \sqrt{2} D(b,b^+)^{1/2} \cos\beta\\
&\leq (1+\epsilon)^3  \tau' 2 \widetilde{R}(b^+,d)  D(b,b^+)^{1/2} \cos\beta.
\end{split}
\end{align}
These estimates  use $\cos \beta \geq 0$ throughout, which holds due to part 1). Moreover, 
line four uses (\ref{26_1}) (c),
line five uses (\ref{26_1}) (a),  line six uses (\ref{26_1}) (b), and the last line uses (\ref{limsup_2}).

3)
Now recall that $\bar{a}-b^+$ is a proximal normal to the set $B$ at $b^+$ with regard to the
$\|\cdot\|_{b^+}$-geometry, with $d = (\bar{a}-b^+)/\|\bar{a}-b^+\|_{b^+}$ the corresponding unit proximal normal. 
Since $\beta = \angle_{b^+} (b-b^+,\bar{a}-b^+) < 90^\circ$ according to part 1),
 we can 
choose $R>0$ such that the point $\hat{a} = b^+ +Rd$ on the proximal normal satisfies
$\|\hat{a} - b^+\|_{b^+} = \|\hat{a}-b\|_{b^+} = R$, so that $b,\hat{a},b^+$ form an equilateral triangle with
two angles $\beta$ adjacent to the side $b,b^+$ of length $\|b-b^+\|_{b^+}$, and two sides of equal length $R$ joining $\hat{a}$.
Therefore $\frac{1}{2}\|b-b^+\|_{b^+} = R \cos\beta$ by the perpendicular bisector theorem. 

Now the $\|\cdot\|_{b^+}$-euclidean ball with center $\hat{a}$ and radius $R$  contains the points $b,b^+\in B$
on  its boundary. By definition of the reach of $B$ at $b^+$ with regard to the norm $\|\cdot\|_{b^+}$, and since the ball in question has its center on the
$\|\cdot\|_{b^+}$-normal $b^++\mathbb R_+d$, this means that $R$ must be at least as large as the $\|\cdot\|_{b^+}$-reach
$\widetilde{R}(b^+,d)$ in that direction. We derive
\begin{equation}
\label{s_second}
\|b-b^+\|_{b^+} = 2R \cos\beta \geq 2\widetilde{R}(b^+,d) \cos\beta. 
\end{equation}
Plugging this into (\ref{27_1})
gives
\begin{align*}
\langle \nabla f(a^+)-\nabla f(b^+),b-b^+\rangle &\leq 
(1+\epsilon)^3 \tau'   \|b-b^+\|_{b^+} D(b,b^+)^{1/2}\\
&\leq (1+\epsilon)^4 \tau' \sqrt{2} D(b,b^+)\\
&< (1-\ell) D(b,b^+)
\end{align*}
by the choice of $\epsilon$ and $\tau'$.
Therefore by the cosine theorem for Bregman distances
\begin{align*}
D(b,a^+) &= D(b^+,a^+) + D(b,b^+) - \langle \nabla f(a^+) - \nabla f(b^+),b-b^+\rangle \\
&\geq D(b^+,a^+) + D(b,b^+) - (1-\ell) D(b,b^+)\\& = D(b^+,a^+) + \ell D(b,b^+). 
\end{align*}
\end{proof}

\begin{remark}
We see from the argument in part 3), and also from a similar one in Proposition \ref{prop10}, that the technique  in part 2) of the proof of Proposition \ref{prop9} looks like
a Bregman version of the euclidean perpendicular bisector theorem.
\end{remark}

\subsection{Convexity}
As we had seen in  Proposition \ref{prop_convex}, convexity of $B$, or $\nabla f(A)$,
makes things easier, but surprisingly, convexity of $A$ doesn't seem to help. This discrepancy was also observed in \cite[Thm. 7.3]{shawn}, where the right Bregman Chebyshev 
condition of $A$  was shown to imply  convexity of $\nabla f(A)$, not of $A$.  For short, convexity is not
amenable to duality.  Can we still get something when $A$ is convex? 

\begin{proposition}
\label{hard}
Suppose $A$ is closed bounded contained in {\rm int(dom$f$)} and has positive reach. Suppose $f$ is of class $C^{2,1}$. Then there exists $r > 0$ such that the $lr$-three-point inequality
is satisfied at all gaps with gap value $r^*\leq r$.
\end{proposition}

\begin{proof}
Since $A$ has positive reach and
$\nabla f$ is a $C^{1,1}$-diffeomorphism, the image $B^*=\nabla f(A)$ has also positive reach.  
Since positive reach is equivalent to positive mobile reach, we have
$\widetilde{R}(b^{*+}) \geq \widetilde{R}$ for some $\widetilde{R}> 0$  for the mobile reach in dual space and all $b^{*+}\in B^*$.
Therefore, applying Proposition \ref{prop10} in dual space, we get the 
$rl$-three-point inequality for dual building blocks $b^* \stackrel{r*}{\longrightarrow} a^{*+} \stackrel{l*}{\longrightarrow} b^{*+}$ near gaps with  
dual gap distance $\rho_* = \|\nabla f^*(b^{*+})-\nabla f^*(a^{*+})\|_{b^{*+}} <\widetilde{R}$.

Going backwards in the dual formula (\ref{dual_formula}), this means every $lr$-building block $a \stackrel{l}{\longrightarrow} b\stackrel{r}{\longrightarrow} a^+$ 
satisfies the $lr$-three-point inequality as soon as 
$\rho_* =  \|\nabla f^*(b^{*+})-\nabla f^*(a^{*+})\|_{b^{*+}} = \|a^+-b\|_{b^{*+}} < \widetilde{R}$.
Now
$\|u\|_{b^{*+}}^2 = \langle \nabla^2f^*(b^{*+})u,u\rangle = \langle\nabla^2f(a^+)^{-1}u,u\rangle$. Hence  in primal space the sufficient condition reads
$\rho_* = \langle \nabla^2f(a^+)^{-1}(a^+-b),a^+-b \rangle^{1/2} < \widetilde{R}$. Since $\nabla^2f(a^+)^{-1} \succeq \epsilon > 0$ on $A$, we have an estimate of the form
$m' D(b,a^+)^{1/2} \leq  \|a^+-b\|_{b^{*+}}$ similar to (\ref{star_bounds}) combined with (\ref{bounds}) with the same $m'$ for $A,B$, so that the $lr$-three-point inequality now holds 
for building blocks with $D(b,a^+)^{1/2} < m'^{-1} \widetilde{R}$, hence it holds for gaps with  $r^* \leq r :=m'^{-1}\widetilde{R}$.
\end{proof}

This shows that it is prox-regularity, or positive reach, which is amenable to duality, not convexity.  However, the distortion caused by $\nabla f$ makes it 
hard to relate the reach of a set $A$ in primal space to the reach of its $\nabla f$-image $B^*$ in dual space. Quantifying $r$ could at best be
achieved in specific situations.

\begin{corollary}
Let $\bar{x} \in A \cap B$ and suppose $A$ is prox-regular at $\bar{x}$. Let $A,B,f$ be definable, and suppose $f$ is of class $C^{2,1}$. Then there exists a neighborhood $V$
of $\bar{x}$ such that every Bregman alternating sequence which enters $V$ converges to a point in $A \cap B$.
\end{corollary}

This should be compared with Corollary \ref{cor2}.

\section{EM algorithm}
\label{sect_em}

We apply our convergence theory to the {EM} algorithm via the variant
in \cite{csiszar1}, known as the {\it em}-algorithm. Criteria for the two to coincide are given in \cite[Thm. 4, Ex. 10, \S 7.2]{amari}.
As observed in \cite{wu}, or \cite[Sect 3.6]{lachlan},
even when convergent, {EM} iterates may go to local minima or saddle points
of the likelihood, which is not surprising in a non-convex setting. 
In this work, we focus on convergence of the iterates, which contains value convergence  first discussed in \cite{wu}, see also \cite{kunstner}. 
As literature on convergence of the {EM} algorithm is vast, we confine ourselves to a few pointers
\cite{csiszar1,csiszar2,amari,amari_book,tseng,brown,byrne1,censor,cen_zak,hino,kunstner,wu}. En excellent survey up to 2008 is \cite{lachlan}.

\subsection{Exposition of the method}
The {\it em}-algorithm requires a complete data space $X$, an incomplete data space $Y$, a measurable mapping
$T:X \to Y$,  a family of probability measures $P$ on $X$ and their images $P^T$ under $T$ on $Y$.
Assuming $P \ll \mu$ for a $\sigma$-finite base measure on $X$, and $P^T \ll \mu^T$, we have
densities $p = \frac{dP}{d\mu}$ and $p^T = \frac{dP^T}{d\mu^T}$. We further require
an empirical distribution $\hat{P}$ on  $Y$ with $\hat{P} \ll \mu^T$, $\hat{p} = \frac{d\hat{P}}{d\mu^T}$,
which we use to
define the data set on $X$ as
$\mathcal D=\{P: P^T = \hat{P}\}$, respectively, $D = \{p: p^T = \hat{p}\}$.
Finally,  we need a statistical model $\mathcal M$ of distributions $Q \ll \mu$ on $X$, $M = \{q: Q\in \mathcal M\}$, and our goal is to minimize the Kullback-Leibler
information  distance \cite{kullback} between $D$ and $M$, given by
$$
K(p||q) = \int_X p(x) \log \frac{p(x)}{q(x)} d\mu (x).
$$ 
Then the {\it em}-algorithm is the following
alternating procedure:
\begin{align*}
\mbox{{\it e}-step} &\qquad p \in \argmin_{p'\in D} K(p'||q) \\
\mbox{{\it m}-step} &\qquad q^+ \in \argmin_{q'\in M} K(p||q')
\end{align*}
The {\it m}-step is similar to the M step of the EM algorithm and represents maximum likelihood estimation in complete data space.
On the other hand,  the {\it e}-step may differ from the E step, as shown in \cite[6.3]{amari}, and according to \cite[Thm. 4]{amari}, both agree iff 
the conditional expectation with respect to a candidate distribution  of the missing data, given the observed data, is an affine function of  observed data.

It turns out that the {\it e}-step is explicit. As follows from \cite[Thm. 4.1]{kull_leib}, we have
$K(p||q) \geq K(p^T||q^T)$, with equality iff $\frac{p(x)}{q(x)} = \frac{p^T(T(x))}{q^T(T(x))}$ $\mu$-a.e. Applying this to the data set $D$ shows that
$\cev{P}_D(q) =p$ is realized by setting $p(x)= {\hat{p}(T(x))}\frac{q(x)}{q^T(T(x))}$ $\mu$-a.e.

\subsection{Discrete measures}
In the first place, let us assume $X=I,Y=J$ finite, with $\mu$ the counting measure. Then
distributions on $I$ are vectors $p=(p_i)_{i \in I}$ with $p_i \geq 0$, $\sum_{i\in I} p_i = 1$, and similarly,
$\hat{p}=(\hat{p}_j)_{j\in J}$ with $\hat{p}_j \geq 0$, $\sum_{j\in J} \hat{p}_j=1$.
The Kullback-Leibler divergence on $\mathbb R^I_+$ 
$$
K(p||q) = \sum_{i\in I} p_i \log \frac{p_i}{q_i}-p_i+q_i,
$$
is now the Bregman divergence generated by  $f(x)=\sum_{i\in I} x_i\log x_i-x_i$, making the algorithm amenable to our convergence theory.
We assume $\hat{p}_j > 0$ for all $j\in J$ and observe that the data set is
$
D = \left\{p\geq 0: \sum_{i\in I} p_i=1,\sum_{T(i)=j} p_i = \hat{p}_j \;\forall j\in J\right\}$. 
\begin{theorem}
Let $M\subset \mathbb R^I_{++}$ be a closed definable set of statistical model distributions. Let
$p^k,q^k$ be a sequence generated by the $em$-algorithm. Then the $p^k$ converge to a distribution $p^*\in D$.
Every accumulation point $q^*\in M$ of the $q^k$ satisfies
\begin{equation}
\label{cond_exp1}
p_i^* =  \hat{p}_j\, \frac{q_i^*}{\sum_{T(i')=j}q_{i'}^*}, \quad i\in T^{-1}(j), j\in J.
\end{equation}
When $M$ is definable in $\mathbb R_{an}$, then the speed of convergence is no worse than
$\|p^k-p^*\| = O(k^{-\rho})$ for some $\rho\in (0,\infty)$.
\end{theorem}

\begin{proof}
We have $p \in \cev{P}_D(q)$ and $q^+\in \vec{P}_M(p)$, hence $M=A$ and $D=B$ in our general scheme. 
Since $D$ is convex and not entirely contained in the boundary of $\mathbb R^I_+$,  it is interiority preserving, hence
can be pre-processed as in Section \ref{interior}. Therefore
$lr$-building blocks
$q\overset{l}{\underset{e}{\longrightarrow}} p \underset{m}{\overset{r}{\longrightarrow}} q^+$ satisfy the three-point inequality
by Proposition \ref{prop_convex}. Moreover,
$D$ is semi-algebraic, $K(p||q)$ is definable because $\log$ is definable in $\mathbb R_{an,\exp}$, and $M$ is definable by hypothesis.
Therefore the  $lr$-angle condition holds, and convergence of the $p^k$ follows with Theorem \ref{theorem3}.

Next observe that the {\it e}-step can be made explicit due to the structure of $D$ and the observation made above. We have
$$
p_i = \hat{p}_j\, \frac{q_i}{\sum_{T(i')=j}q_{i'}}, \quad i\in T^{-1}(j), j\in J
$$
for $q \stackrel{l}{\longrightarrow} p$, and
we may clearly pass to the limit in this expression, which gives (\ref{cond_exp1}).

Finally, assume $M$ is definable in $\mathbb R_{an}$. Since
$D$ is semi-algebraic, hence  definable in $\mathbb R_{an}$, it remains to show definability of $K(p||q)$ in $\mathbb R_{an}$.
Now
observe that
by $M \subset \mathbb (0,1]^I$ and closedness of $M$ the $q_i^k$ stay bounded away from $0$, i.e. $q_i^k \geq \epsilon > 0$ for all $i,k$. Since 
$\hat{p}_j > 0$ for all $j\in J$, formula (\ref{cond_exp1}) shows that all $p_i^k$ also stay uniformly bounded away from $0$, say $p_i^k \geq \delta > 0$
for all $i,k$. Bearing in mind that $p_i^k \leq 1$,
the logarithms $\log p_i$ in $K(p||q)$ can a priori be restricted to the  interval $[\delta,1]$. But the restriction $\log\upharpoonright [\delta,1]$ is globally sub-analytic, 
cf. \cite{dries,dries_miller}.
In consequence
the K\L-inequality becomes a \L ojasiewicz inequality with de-singularizing function
$\phi'(s) = s^{-\theta}$ for some $\theta\in [\frac{1}{2},1)$. Therefore Corollary \ref{cor_speed}  gives the claimed rate of convergence.
\end{proof}

\begin{remark}
In the feasible case $D\cap M$ has a neighborhood of attraction $W$ such that any sequence
$p^k,q^k$ entering $W$ converges to some $p^* \in D \cap M$ with the same rate
$\|q^k-p^*\| = O(k^{-\rho})$, $\|p^k-p^*\|=O(k^{-\rho})$. Recall from Corollary \ref{cor_speed} that $\rho$ can be expressed in terms of the Lyapunov exponent $\theta$
of $F$, hence in some cases an even more specific rate may be obtainable.

Convergence of the $q^k$ with the same rate also occurs for gaps $>0$
when the left Bregman reach of $A=M$ at the $q^k$ is larger than the gap value, as $\cev{P}_M$ is then locally Lipschitz.  
\end{remark}

Note that (\ref{cond_exp1}) says  $p^k = \mathbb E_{q^k}(p|\hat{p})$, i.e., $p^k$ is the conditional expectation of $p$ given $\hat{p}$ with regard to the probability distribution
$q^k$. 
In the limit we then have $p^*=\mathbb E_{q^*}(p|\hat{p})$ for all accumulation points $q^*$ of the $q^k$.
 Since definability is not a severe restriction, this
is quite useful in practice.

\subsection{Exponential family}
We consider
an exponential  family of densities with respect to a base measure $dx$
\begin{equation}
\label{exp_family}
p_\theta(x) = \exp( \langle \theta, t(x)\rangle - f(\theta) + k(x)),
\end{equation}
where $t(x)$ is the sufficient statistic, $\theta$ the natural parameter varying in a parameter set $\Theta$, $f(\theta)$ the log-normalizer function, 
and $\exp k(x)$ the carrier measure density. 
\begin{lemma}
\label{Kull_Breg}
The Kullback-Leibler divergence of two distributions $p_\theta(x)$ and $p_{\theta'}(x)$ belonging to the same exponential family is
$K(p_{\theta'}||p_{\theta}) = D_f(\theta,\theta')$, where $D_f$ is the Bregman divergence induced by the log-normalizer $f$.
\end{lemma}

\begin{proof}
From (\ref{exp_family}), and since $\int_Xp_\theta(x) dx=1$,  we have
$$
f(\theta) = \log \int_{X} \exp\{ \langle \theta,t(x)\rangle + k(x)  \} dx.
$$
Differentiation with respect to $\theta$ gives
$$
\nabla f(\theta) = {\int_X t(x) \exp\{ \langle \theta,t(x)\rangle+k(x)\} dx }\bigg/{\int_{X} \exp\{ \langle \theta,t(x)\rangle+k(x) \}dx}.
$$
Now $\exp f(\theta) = \int_X \exp\{  \langle \theta, t(x)\rangle + k(x) \} dx$, hence $\nabla f(\theta) =\int_X t(x) \exp\{ \langle \theta, t(x)\rangle - f(\theta) + k(x)  \} dx
= \int_X t(x) p_\theta(x) dx=\mathbb E_\theta [t(x)]$, the expectation  of the random variable $t(x)$ with respect to the distribution
$p_\theta dx$ (see also \cite[(1.57)]{lachlan}). Then
\begin{align}
\label{reversed}
\begin{split}
K(p_{\theta'}||p_{\theta}) &= \int p_{\theta'}(x) \log\frac{p_{\theta'}(x)}{p_{\theta}(x)}\, dx \\
&= \int_X p_{\theta'}(x) \left( f(\theta)-f(\theta') + \langle \theta'-\theta,t(x)\rangle  \right) dx\\
&=\int_X p_{\theta'}(x) \left(  D_f(\theta,\theta') + \langle \theta-\theta',\nabla f(\theta')\rangle + \langle \theta'-\theta,t(x)\rangle \right) dx\\
&=D_f(\theta,\theta') + \int_X p_{\theta'}(x) \langle \theta-\theta', \nabla f(\theta')-t(x)\rangle dx\\
&= D_f(\theta,\theta') + \langle \theta-\theta', \nabla f(\theta') - \mathbb E_{\theta'}[t(x)]\rangle \\
&= D_f(\theta,\theta').
\end{split}
\end{align}
\end{proof}
In order to connect with our general  set-up, we have to make sure that $f$ is Legendre. We have the 

\begin{definition}
{\rm The exponential family is called {\it steep} if its log-normalizer $f$ is of Legendre type. }
\end{definition}  

An exponential family is {\it regular} if the natural parameter space $\Theta$ is open.
It is known that regular exponential families are steep, but the steep class is larger; cf. \cite[Ex. 3.4]{brown}, \cite{banerjee}.
We are now ready to concretize the $em$-algorithm for exponential families,
specifying 
model set
$\mathcal M = \{p_\theta: \theta \in M\}$ and data set $\mathcal D=\{p_\theta: \theta \in D\}$
by their parameter representatives $M,D\subset \Theta$.

\begin{algorithm}[h!]
\caption{$em$-algorithm for exponential family}\label{algo_em_2}
\noindent\fbox{%
\begin{minipage}[b]{\dimexpr\textwidth-\algorithmicindent\relax}
\begin{algorithmic}
\STEP{$e$-step} Given current model density $p_{\theta}$, $\theta\in M$, complete data with the help of the sample via $\theta'\in \vec{P}_D(\theta)$.
Obtain completed data density $p_{\theta'}$, $\theta'\in D$. 
\STEP{$m$-step} Given complete data  density $p_{\theta'}$, $\theta'\in D$, improve parameter estimate by maximum likelihood in complete data space via
$\theta \in \cev{P}_M(\theta')$. Back to step 1.
	\end{algorithmic}
\end{minipage}
}%
\end{algorithm}

We observe that due to (\ref{reversed}) left and right projections have been swapped, and
the algorithm has now the form $\theta'\underset{m}{\overset{l}{\longrightarrow}}\theta \underset{e}{\overset{r}{\longrightarrow}} \theta^{'+}\underset{m}{\overset{l}{\longrightarrow}}\theta^+$,
which matches (\ref{index_free}) with
$M=B$, $D=A$, $a=\theta'$, $b=\theta$, $a^+=\theta'^+$, $b^+= \theta^+$.

It remains to say a bit more about the data set. Following  \cite{amari}, we consider 
the expectation parameter $\eta(\theta) = \mathbb  E_\theta[t(x)]$, which by Lemma \ref{Kull_Breg} is 
$\eta = \nabla f(\theta)$, with inverse $\theta = \nabla f^*(\eta)$. When
the sufficient statistic is of the form $t(x) = (t_1(x),t_2(x))$ for observed $y=t_1(x)$ and hidden  $z=t_2(x)$,  then 
\begin{equation}
\label{affine}
D=\{\theta\in G: \mathbb E_\theta[t_1(x)] = \hat{y}\}, \quad \mathcal D = \{p_\theta: \theta\in D\},  
\end{equation}
where $\hat{y}$ is the available sample. We can partition $\theta = (\theta_1,\theta_2)$ accordingly, so that
$\langle \theta,\eta\rangle = \langle \theta_1,\eta_1\rangle + \langle \theta_2,\eta_2\rangle$, then 
$\nabla f(D)=\{(\eta_1,\eta_2)\in G^*: \eta_1 = \hat{y}, \eta_2 \mbox{ free}\}$, which in expectation coordinates is an affine
subspace $L^*$, intersected with dom$\nabla f^*=G^*$. 
\begin{theorem}
\label{thm_6}
Suppose the log-normalizer  $f(\theta)$ of the steep exponential family 
and the model parameter set
$M$ are definable. Suppose $M\subset G$ is closed bounded and
$D$ has the structure {\rm (\ref{affine})}. Then the 
$em$-algorithm $\theta'^+ = \vec{P}_D \circ \cev{P}_M(\theta')$ converges. 
\end{theorem}

\begin{proof}
1)
By hypothesis $M$ is closed bounded and contained in $G$, so its image $A^*:=\nabla f(M)$ is closed bounded and contained in $G^*$. 
Since $D \subset G$, we have $B^*:=\nabla f(D) \subset G^*$, and by the above $B^*=L^*\cap G^*$.
This means $B^*=\nabla f(D)$ is convex, while not necessarily closed.

2)
Put $B_1^*=L^* \cap {\rm cl}\,G^*$, then $B_1^*$ is closed and satisfies the constraint qualification $B_1^* \cap G^* = B^*\not=\emptyset$. Therefore
$B_1^*$ is interiority preserving (Section \ref{interior}, \cite[Thm. 3.12]{BB_legendre}), hence $\mbox{${\cev{P}}^{{}^{{}^*}}_{B_1^*}(y^*)$}$ is defined for $y^*\in G^*$
and we have $\mbox{${\cev{P}}^{{}^{{}^*}}_{B_1^*}(y^*)$}\subset B_1^*\cap G^*=B^*$. This little detour shows that
$\mbox{${\cev{P}}^{{}^{{}^*}}_{B^*}(y^*)$}\not=\emptyset$ is well-defined for $y^*\in G^*$.

3)
By compactness of $A^*\subset G^*$ dual right projections on $A^*$ are well defined, and by 2) dual left projections from $A^*$ are also well-defined and, moreover,
go to $B^*$. Now consider $\theta'\in D$, then $\theta \in \cev{P}_M(\theta')$ is well defined by compactness of $M\subset G$, i.e. $\theta' \stackrel{l}{\longrightarrow} \theta$
is well defined in primal space. Hence $\nabla f(\theta') \stackrel{r*}{\longrightarrow} \nabla f(\theta)$ is well defined in $G^*$. But by what we had just seen we can now continue left
projecting from $\nabla f(\theta)\in A^*$  into   $B^* \subset G^*$, hence the dual $rl$-building block
$\nabla f(\theta') \stackrel{r*}{\longrightarrow} \nabla f(\theta)\stackrel{l*}{\longrightarrow}\nabla f(\theta'^+)$ is well defined and lies in $G^*$. By duality backward,
the primal $lr$-building block
$\theta' \stackrel{l}{\longrightarrow} \theta \stackrel{r}{\longrightarrow} \theta'^+$ is now also well defined and lies in $G$. 
Iterating this, the entire primal sequence is well-defined and lies in $G$, and its mirror sequence lies in $G^*$. It remains to show that
the same holds for the accumulation points of these sequences.

4)
As $B_1^*$ is interiority preserving,
we may apply Section \ref{interior} to the dual alternating sequence, which means there exists a closed bounded subset $C^*$ of $B_1^*\cap G^*=B^*$ such that 
left projections of iterates from $A^*$ go to $C^*$. In consequence the dual sequence alternates between $C^*$ and $A^*$, 
now including accumulation points. Mapping this back via  $\nabla f^*$
means the primal sequence alternates between $M=\nabla f^*(A^*)$ and the compact $C:=\nabla f^*(C^*) \subset D$, and accumulation points belong to $M$, respectively, $C$. 
This makes the situation amenable to our convergence theory.

5)
For that it remains to establish the angle condition.
Since $f,M$ are definable by hypothesis, it remains to check
definability of $C$. Now observe that $G$ as the interior of the domain of $f$ is definable.  Definability of $\nabla f$ also follows from definability of $f$.
Therefore $G^* = \nabla f(G)$ as the image of a definable set under a definable diffeomorphism is definable (see \cite[Prop. 1.6]{coste}). Hence
$B^* = L^* \cap G^*$ is definable, $L^*$ being algebraic. Now recall that the construction in Section \ref{interior} gives a definable  $C^*$ because $B^*$ is definable, and hence
$C = \nabla f^*(C^*)$ is also definable, using that $\nabla f^*=(\nabla f)^{-1}$ is also a definable diffeomorphism. This means
the $lr$-angle condition is satisfied. 

Now all the hypotheses of  Theorem \ref{theorem3} are satisfied, which gives convergence of the primal $lr$-sequence $\theta'^+ = \vec{P}_D \circ \cev{P}_M(\theta')$.
\end{proof}

\begin{remark}
1) The argument hinges on $M$ being bounded, which is not always true in practice, but some boundedness hypothesis is required
(see e.g. \cite[(6)]{wu}), because in the infeasible case even euclidean alternating projections between unbounded convex sets
may escape to infinity, and without convexity, this may happen even in the feasible case. 
When the sequence $\theta_k,\theta'_k$ is bounded and the $\theta'_k$ stay away from $\partial G$, we may always select a closed bounded subset $M_0\subset G$ of the model set $M$
such that $a_k,b_k$ remain alternating between $D,M_0$, and then in a second step use the treatment of Section \ref{interior} to get a bounded
$D_0\subset G$.

2) The result holds more generally when $y=t_1(x)$ is affine in $z=t_2(x)$, $y = Az + b$. Let $A \in \mathbb R^{n\times m}$ of maximal rank $m < n$,
find $Q$ invertible $n\times n$ with $A = [\widetilde{A} \; 0]\, Q$ and $\widetilde{A}$ regular of size $m \times m$, and
make the change of coordinates $\tilde{z} = {\rm diag}(\widetilde{A} , I_{n-m}) Q z =: Tz$, then $\tilde{z}=(\tilde{y},\tilde{v})$ with $\tilde{y}=\pi_1(\tilde{z}) = Az$ and $y=\tilde{y}+b$,
which reduces the affine case to (\ref{affine}). 

3) 
An intriguing question is whether  EM, respectively {\it em}, may fail to converge and generate a continuum of accumulation points. 
A natural place to look are euclidean alternating projections, as those arise when
estimating the mean of a Gaussian with known variance. Counterexamples for AP with a continuum of accumulation points have been given in \cite{bauschke,DR}, 
but 
do not apply to EM, because in these examples the structure of the set $D$ playing the role of the  data set is too exotic. We sketch a possible counterexample in Section \ref{sect_examples}.

4)
Inspecting lists of exponential families, one finds that log-normalizers $f$ often feature terms  $\log \theta_i$ for components of the natural parameter $\theta$. 
All cases we are aware of are governed by the  o-minimal structure $\mathbb R_{an,\exp}$ unifying globally sub-analytic sets with exponential and logarithm;
cf. \cite{dries,dries_book,wilkie}. Yet there is interest to arrange model parameters $\theta\in M$ such that in these
$\log \theta_i$-terms the $\theta_i$ may {\it a priori} be bounded on some $[\underline{\theta},\overline{\theta}]$ with $0 < \underline{\theta}$, $\overline{\theta} < \infty$. Namely,
by \cite{dries,dries_miller}, this has the benefit that
$f$ will  be definable in  $\mathbb R_{an}$. Since $D$ is by default definable in $\mathbb R_{an}$, things depend on $M$.
When $M$ is also definable in $\mathbb R_{an}$, we get
de-singularizing functions
$\phi(s) = s^\alpha$ for some $\alpha \in [\frac{1}{2},1)$, which by
Corollary \ref{cor_speed} allows to quantify the speed of convergence to some $O(k^{-\rho})$.

5) 
Instances of sublinear speed of EM are mentioned e.g. in \cite[p. 102]{lachlan}, even though the general understanding seems to be that EM should converge linearly.
However, it may be extremely hard to predict linear convergence a priori. For instance,
even in the feasible case $p^*\in D \cap M$, and despite the convenient structure of $D$, we would have to show that $D,M$ intersect transversally at $p^*$ prior to coming to know $p^*$.
This is possible only in very specific situations. Realistically, we should therefore only claim a rate $O(k^{-\rho})$.  That may be predicted credibly, as
definability of $D,M,f$ is usually easy to check. A fair chance to prove transversality might be the case (\ref{affine}) when $M$
has a simple structure.
\end{remark}

\begin{remark}
In those cases where complete data are $x=(y,z)$ with observed $y$ and hidden $z$,
the data set is chosen as
\begin{equation}
\label{new}
D= \{\theta \in G: \eta = t(\hat{y},z), \mbox{$\hat{y}$ observed, $z$ arbitrary}, \theta = \nabla f^*(\eta)\}.
\end{equation}
Here the {\it em}-algorithm converges under the assumptions of Theorem \ref{thm_6} if  $D$ now with structure (\ref{new})  is definable,
which it is e.g. if  the mapping $t(\hat{y},\cdot)$ is definable (cf. \cite[Prop. 1.6]{coste}).

The following modification
gives the corresponding result for the EM algorithm. Here we have to re-define the data set as
\begin{equation}
\label{newnew}
D = \{\theta\in G: \eta = \mathbb E\left( t({y},z) \big | \hat{y}, q(z|\hat{y}) \right), \theta = \nabla f^*(\eta)\},
\end{equation}
where $q(z|y)$ is an arbitrary distribution of the imputed value $z$, given $y$, and the conditional expectation is with regard to the measure
$q(z|\hat{y}) dz$. Then the EM algorithm is the method of alternating Bregman projections
between $M$ and this modified $D$, and Theorem \ref{thm_6}
assures convergence as soon as $D$ has structure (\ref{newnew}) and  is definable. Here definability of $t$ alone is not sufficient, unless $t(\hat{y},\cdot)$ is affine in $z$. In general a case-by-case analysis
seems necessary.
\end{remark}

\subsection{dSPECT imaging}
We end with an application in dynamic SPECT (dSPECT) imaging \cite{dynEM}. Voxels $i\in I$ have unknown activity $x_{ik}$ varying in time
$t_k$,  $k\in K$.
Camera bin $j\in J$  receives  $y_{jk}$ counts at time $t_k$ at angular position $\alpha_k$,  where $\mathbb E(y_{jk}) = \sum_{i\in I} c_{ijk} x_{ik}$, and the known coefficients 
$c_{ijk}$ reflect camera and collimator geometry. Complete data $z$ are activities
$z_{ijk}$ emitted from voxel $i$ at time $t_k$ to camera bin $j$ in position $\alpha_k$, $\mathbb E(z_{ijk}) = c_{ijk} x_{ik}$. The $x_{ik}$ are the unknown parameters. 
Two laws have been discussed in \cite{dynEM}, Poisson, and Gaussians with known variance. 
The data set is $D = \{z: \sum_{i\in I} z_{ijk} = y_{jk} \forall j\in J, k\in K\}$. Since the problem is underdetermined, prior information on the time behavior 
of the $x_{ik}$ is added. In \cite{dynEM} the authors use $x_{ik} = a_i e^{-\lambda_i t_k} + b_i e^{-\mu_it_k} + d_i$, and this defines a model
$M_e=\{v: v_{ijk} = c_{ijk} x_{ik}(a_i,b_i,d_i,\lambda_i,\mu_i)\}$ with $5|I|$ parameters. The variant in \cite{maeght} uses a Prony model 
$x_{ik} = \alpha_{i} x_{i,k-2} + \beta_{i} x_{i,k-1} + \gamma_{i}$ with $M_p=\{v: v_{ijk} = c_{ijk} x_{ik}(\alpha_i,\beta_i,\gamma_i,x_i^0,x_i^1)\}$, where the $5|I|$ parameters are
$\alpha_i,\beta_i,\gamma_i$ and two initial values $x_i^0$, $x_i^1$ per voxel. Clearly $M_p$ is definable in $\mathbb R_{an}$, 
$M_e$ in $\mathbb R_{an,\exp}$, so that convergence of the E step sequence for both cases is assured,
giving convergence of the corresponding $x$. In the feasible case the M sequence converges as well. For the Prony model $M_p$
the rate is $O(k^{-\rho})$, for  $M_e$ this holds when the exponentials can be restricted to a compact interval. 
For the Gaussian case convergence could be obtained from \cite{gerchberg,noll}, while for the Poisson model
convergence is now for the first time established here.

\section{Examples}
\label{sect_examples}

\begin{example}
\label{ex_2}
We consider the Kullback-Leibler divergence in $\mathbb R^2$,
$K(x||y) = \sum_{i=1}^2 x_i \ln(x_i/y_i) -x_i + y_i$, with $0\ln 0 = 0$.
Take $\bar{y}=(1,1)$, $\bar{x}=(1,0)$ then $K(\bar{x}||\bar{y})=K((1,0)||(1,1)) = 1 = \frac{1}{2} r^2$ with $r=\sqrt{2}$. The
Bregman ball $\cev{\mathcal B}(\bar{y},\sqrt{2})$ contains $(1,0)$, but no other point $(z,0)$, $z \not=1$, because
$K((z,0)||(1,1)) = z\log z - z+2 \stackrel{!}{=} 1$ has only the solution $z=1$, while $z\log z - z + 2 > 1$ for $z\not=1$. 
Hence the line $x_2=0$ is tangent (a support hyperplane) to $\cev{\mathcal B}(\bar{y},\sqrt{2})$
at $(1,0)$.

We now squeeze a curve $B$ in between
$x_2=0$ and $\partial \cev{\mathcal B}(\bar{y},\sqrt{2})$ in such a way that $B \cap \{x_2=0\} = \{\bar{x}\} = B \cap \cev{\mathcal B}(\bar{y},\sqrt{2})$. Then $B\cap G \not=\emptyset$, but
$\cev{P}_B(\bar{y}) = \bar{x} \not\in  G=\mathbb R^2_{++}$. Make the ansatz $B = \{(z,f(z)):  1-\epsilon \leq z \leq 1+\epsilon\}$ with $f(1)=1$, then
$K((z,f(z))||(1,1)) = z\log z - z +  f(z) \log f(z)-f(z) + 2 > 1$ when we arrange $f(z)$ such that
$f(z) \log f(z)-f(z) \leq \frac{1}{2} \left[ -1 -z\log z + z  \right]$ on $(1-\epsilon,1+\epsilon)$. 
\end{example}

\begin{example}
Failure of convergence of euclidean alternating projections with at least one of $A,B$ non-convex are 
given e.g. in \cite{bauschke,DR,noll,gerchberg}. In those cases $a_k,b_k$ are bounded with $a_{k-1}-a_k\to 0$, $b_{k-1}-b_k \to 0$, but fail to converge because either
angle condition or regularity fail, while the other holds, producing a continuum of accumulation points. For pictures see \cite{bauschke,DR}.
\end{example}

\begin{example}
We sketch an example of failure of convergence of a bounded alternating sequence with a continuum of accumulation points in the feasible case $A \cap B\not=\emptyset$,
where one of the sets is affine. This could be re-organized to give EM
for estimating the mean of a gaussian with known variance with  a non-convex curved parameter set. 

Let $A=\{(x,0): x\in \mathbb R^n\}$ and $B= \{(x,f(x)): x\in \mathbb R^n\}$ the graph of a $C^1$-function
$f:\mathbb R^n\to \mathbb R$ with $f(x) >0$ on $|x| < 1$, $f(x)=0$ on $|x|=1$. 
Take euclidean alternating projections,  then $P_A(x_k,f(x_k)) = (x_k,0)$, while
$(x_k,f(x_k)) \in P_B(x_{k-1},0)$ iff  $x_{k-1} = x_k + f(x_k) \nabla f(x_k)$. This method follows steepest ascent
backwards. 

Taking $n=2$, 
we let $B$  the graph of the mexican hat function \cite{absil} on $x_1^2+x_2^2 \leq 1$, or likewise, its epigraph.
Similar to the argument given for steepest descent with infinitesimal steps in \cite{absil},  AP with infinitesimal steps will also follow the valley of the hat downward, 
endlessly circling around and approaching the boundary curve $x_1^2+x_2^2=1$, where 
$f=0$.  Since the stepsize $f(x_k)$ goes in fact to 0, this argument is plausible. 
What is amiss for convergence is the angle condition, the graph of the mexican hat  failing  K\L, while regularity is
guaranteed since $A$ is a plane. For a picture see \cite{absil}.

\end{example}

\begin{example}
\label{slow_reach}
We explain vanishing reach in the euclidean case. 
Let $B=\{(x,|x|^{3/2}): x\in \mathbb R\}$,
then $B$ has vanishing reach at the origin in direction $d=(0,1)$. As shown in \cite{gerchberg}, the radius $R_x$ of the largest ball
touching $B$ at $b=(x,|x|^{3/2})$ from above is of the order $R_x=O(|x|^{1/2})$ as $x\to 0$. In particular, the point $(0,0)\in B$ cannot be projected on from above, while all
other $(x,|x|^{3/2})$ can.
For more details on this example see \cite{noll}.
\end{example}

\begin{example}
In \cite[Ex. 10]{amari} the author presents the case of a curved exponential family where EM and {\it em}-algorithms converge to different limit points, even though these agree asymptotically for large sample sizes. 
The {\it em} algorithm converges to a point in $D \cap M$, while EM converges to a local minimum.
\end{example}

\end{document}